\documentclass[onefignum,onetabnum]{siamart190516}

\newcommand{\RR}{\mathbb{R}}
\newcommand{\CC}{\mathbb{C}}
\newcommand{\rank}{\textnormal{rank}}
\newcommand{\nullspace}{\textnormal{null}}

\usepackage{paralist}
\usepackage{drawmatrix}
\usepackage{tikz}
\usetikzlibrary{shapes,snakes}
\usetikzlibrary{matrix,positioning,decorations.pathreplacing}
\usetikzlibrary{intersections}
\usetikzlibrary{calc}
\usetikzlibrary{arrows}
\usetikzlibrary{fit}

\usepackage{amsfonts}
\usepackage{graphicx}
\usepackage{epstopdf}
\usepackage{algorithmic}
\ifpdf
  \DeclareGraphicsExtensions{.eps,.pdf,.png,.jpg}
\else
  \DeclareGraphicsExtensions{.eps}
\fi

\newsiamremark{example}{Example}
\newsiamremark{remark}{Remark}
\newsiamremark{hypothesis}{Hypothesis}
\crefname{hypothesis}{Hypothesis}{Hypotheses}
\newsiamthm{claim}{Claim}

\headers{A rational Even-IRA algorithm}{P. Benner, H. Fassbender, and P. Saltenberger}

\title{A rational Even-IRA algorithm for the solution of $T$-even polynomial
eigenvalue problems\thanks{Version of September 12, 2020.}}

\author{Peter Benner\thanks{Max-Planck-Institute for Dynamics of Complex Technical Systems Magdeburg, (\email{benner@mpi-magdeburg.mpg.de}.)}
\and Heike Fassbender\thanks{Institute for Numerical Analysis, TU Braunschweig
  (\email{h.fassbender@tu-braunschweig.de}, \email{philip.saltenberger@tu-braunschweig.de}.)}
\and Philip Saltenberger\footnotemark[3]}
\usepackage{amsopn}

\ifpdf
\hypersetup{
  pdftitle={A rational Even-IRA algorithm for the solution of $T$-even polynomial
eigenvalue problems},
  pdfauthor={P. Saltenberger}
}
\fi

\begin{document}

\maketitle

\begin{abstract}
In this work we present a rational Krylov subspace method for solving real large-scale polynomial eigenvalue problems with $T$-even (that is, symmetric/skew-symmetric) structure. Our method is based on the \textsc{Even-IRA} algorithm \cite{EvenIRA}. To preserve the structure, a sparse $T$-even linearization from the class of block minimal bases pencils is applied, see \cite{Dop18}. Due to this linearization, the Krylov basis vectors can be computed in a cheap way. Based on the ideas developed in \cite{BennEff}, a rational decomposition is derived so that our method explicitly allows for changes of the shift during the iteration. This leads to a method that is able to compute parts of the spectrum of a $T$-even matrix polynomial in a fast and reliable way.
\end{abstract}



\section{Introduction}
\label{sec:intro}

Eigenvalue problems are ubiquitous in engineering, physics, mechanics and many more scientific disciplines. Moreover, they lie at the heart of numerical linear algebra.
As eigenproblems stemming from real-world-app\-li\-ca\-tions are often subject to physical constraints and side conditions, they frequently and naturally inherit structure. For instance, mechanical vibration systems are usually described by symmetric mass, damping and stiffness matrices, see \cite{MeerTiss}. Optimal control problems often involve Hamiltonian/skew-Hamiltonian matrix pencils \cite{MehrWat01}.
But, after all, which features and properties single out faithful numerical algorithms for structured problems from universal methods?
In the first place, the occurrence of structure can be utilized to speed up algorithms and reduce memory requirements. This originates from the deeper focus on the true nature of the problem compared to standard methods. In addition to that, the adequate exploitation of structure is beneficial (and indispensable, sometimes) for the reliability of an algorithm. Indeed, a proper numerical treatment of structure will often produce more accurate and physically meaningful results. Consequently, it seems reasonable to design tailor-made algorithms instead of addressing structured problems without any care by standard means.
We present an algorithm for real, $T$-even polynomial eigenvalue problems of large scale that takes into account all the aforementioned aspects.
The method we propose is an implicitly-restarted rational Krylov-Schur approach based on the \textsc{Even-IRA} algorithm introduced in  \cite{EvenIRA} (see also \cite{FaSa18}). In contrast to the \textsc{Even-IRA} algorithm and motivated by \cite{BennEff}, our approach explicitly allows for changes of the shift parameter during the iteration. This leads to a flexible and adjustable rational Krylov algorithm.

There exist various major applications, including the vibration of gyroscopic systems and control theory, that lead to $T$-even polynomial eigenproblems of large size, see e.g. \cite{EvenIRA, Bet13} and the references therein.  A \emph{matrix polynomial} $P(\lambda)$ is an element from $\mathbb{R}^{m \times n}[\lambda]$, i.e.
\begin{equation}
P(\lambda) =  \sum_{k=0}^\ell P_k \lambda^k = P_\ell \lambda^\ell + P_{\ell - 1} \lambda^{\ell -1} + \cdots + P_1 \lambda + P_0 \label{equ:MatPol}
\end{equation}
with matrices $P_j \in \mathbb{R}^{m \times n}$. The degree $\textnormal{deg}(P)$ of $P(\lambda)$ is the largest index $j$ with $P_j \neq 0$. Often, we write $P(\lambda)$ as a matrix with polynomial entries, i.e., as an element from $\mathbb{R}[\lambda]^{m \times n}$. Here, we are mostly interested in square matrix polynomials $P(\lambda) \in \mathbb{R}[\lambda]^{n \times n}$ with some particular structure in its matrix coefficients. We call $P(\lambda) \in \mathbb{R}^{n \times n}[\lambda]$ as in \eqref{equ:MatPol} $T$-even if $P_j = P_j^T$ holds whenever $j$ is even and $P_j = -P_j^T$ holds otherwise. Equivalently, $P(\lambda)^T= P(-\lambda)$. Eigenvalue/eigenvector pairs $(\mu,x) \in \mathbb{C} \times \mathbb{C}^n$ of $P(\lambda)$ are characterized by the relation $P(\mu)x=0$. 
To find eigenvalues of $P(\lambda)$, it is a common approach to turn $P(\lambda)$ into a matrix polynomial $\mathcal{L}_P(\lambda) = \lambda X+Y$ of degree one\footnote{Matrix polynomials of degree one are often called matrix pencils.}  (e.g. the Frobenius companion form, \cite{MMMM06}) by linearization. Then, the eigenvalues of $P(\lambda)$ and $\mathcal{L}_P(\lambda)$ coincide and the generalized eigenproblem corresponding to the linearization $\mathcal{L}_P(\lambda)$ may be solved by, e.g., the standard QZ algorithm, cf. \cite{MolStew73}.
However, solving a structured (i.e. $T$-even) eigenvalue problem via the QZ algorithm and the Frobenius companion form is not conducive in the light of the problems nature and structure.

In particular, the spectrum of $T$-even matrix polynomials has a \emph{Hamiltonian structure}, that is, it is symmetric with respect to both the real and the imaginary axis.
The algorithm we present takes care of this fact in two different ways. On the one hand, the linearization $\mathcal{L}_P(\lambda)=\lambda X + Y$ of $P(\lambda)$ we consider is a symmetric/skew-symmetric matrix pencil (i.e. $Y=Y^T, X= -X^T$). In particular, $\mathcal{L}_P(\lambda)$ itself is $T$-even and
so it naturally preserves the Hamiltonian spectral structure of $P(\lambda)$.  On the other hand, for any $\zeta \in \mathbb{C}$ outside the spectrum of $P(\lambda)$, we consider the special shift-and-invert transformation
$$ \mathcal{L}_P(\zeta) = \zeta X + Y \; \mapsto \; K(\zeta) := \mathcal{L}_P(\zeta)^{-T}X \mathcal{L}_P(\zeta)^{-1}X$$
 as proposed in \cite{MehrWat01, EvenIRA}. Each eigenvalue pair $(+ \mu, - \mu)$ of $\mathcal{L}_P(\lambda)$ is transformed to only one eigenvalue  $\theta = (\mu^2 - \zeta^2)^{-1}$ of $K(\zeta)$. Consequently, $K(\zeta)$ preserves eigenvalue pairings and the spectral symmetry inherent to the problem is respected.

The foundation of our method is the \textsc{Even-IRA} algorithm from  \cite{EvenIRA}. This method is a sophisticated variant of the Krylov-Schur algorithm (see Stewart \cite{KrylovSchur}) applied to $K(\zeta)$ for some appropriately chosen shift parameter $\zeta$ and a $T$-even linearization $\mathcal{L}_P(\lambda)$ for $P(\lambda)$. Rather than applying a (structure-preserving) symplectic Lanczos process as in \cite{BenFaStoll, BenFS11}, our approach is related to the ideas established for the SHIRA algorithm in \cite{MehrWat01} (see also \cite{BennEff}). To define $K(\zeta)$, we take $\mathcal{L}_P(\lambda) = \lambda X + Y$ to be a special linearization from the class of block minimal bases pencils, see \cite{Dop18}. Due to the structure and sparsity of $\mathcal{L}_P(\lambda)$, the computation of matrix-vector-products $K(\zeta)x$ can be realized implicitly without ever forming $K(\zeta)$ at all. Moreover, linear systems with $\mathcal{L}_P(\zeta)$ and $\mathcal{L}_P(\zeta)^T$ (that arise in Arnoldi-like processes from matrix-vector-products involving $K(\zeta)$) can be solved implicitly through systems  involving only $P(\zeta)$ and $P(-\zeta)$. Accordingly, the complexity of computing  $K(\zeta)x$ is reduced by a significant amount since the size of $P(\zeta)$ is substantially smaller than the size of $\mathcal{L}_P(\zeta)$. For the same reason, memory requirements (e.g. for storing matrix decompositions) can be decreased. These advantages of $\mathcal{L}_P(\lambda)$ over other linearizations (see, e.g., \cite{Mack_Diss}) have already been successfully applied in \cite{FaSa18} to the \textsc{Even-IRA} algorithm. However, as the \textsc{Even-IRA} algorithm does not allow for changes of the shift $\zeta$ during the iteration, this feature is incorporated in our method. Based on \cite{BennEff} and \cite{Ruhe84}, our rational \textsc{Even-IRA} algorithm permits shift adjustments during the iteration without discarding the information that has been accumulated so far. Retaining the advantageous computational aspects, this endows our approach with more flexibility.
In consequence, the rational \textsc{Even-IRA} algorithm we present is a new reliable, flexible and fast numerical method with reasonable costs.

This work is structured as follows:
\begin{compactenum}
\item The basic definitions regarding matrix polynomials and their eigenvalues are presented in Section \ref{sec:Notation}. We introduce the concept of linearization and show how a $T$-even linearization can be constructed.
\item In Section \ref{sec:EvenIRA}, we briefly review the \textsc{Even-IRA} algorithm from \cite{EvenIRA}. It is the basis of our rational method to compute eigenvalues of $T$-even matrix polynomials in a structure-preserving way.
\item We show how the matrix-vector-multiplications involved in the \textsc{Even-IRA} algorithm can be carried out in a very efficient and implicit way in Section \ref{sec:LinearSystems}. This is possible without forming the corresponding large-scale matrix at all.
\item Section \ref{sec:RatArnoldiDecomp} is dedicated to the rational Arnoldi decomposition. We show how a rational decomposition can be invoked for the \textsc{Even-IRA} algorithm and how it is applied in a useful fashion for our purpose.
\item We introduce the rational \textsc{Even-IRA} algorithm in Section \ref{sec:RatEvenIRA}. We discuss the Krylov-Schur-restart procedure in detail and also address the issue of infinite eigenvalues to guarantee a stable convergence of the algorithm.
\item Some numerical examples are given in Section \ref{sec:examples}. We also illustrate how the shift-strategy influences the algorithms success.
\end{compactenum}
Some conclusions are given in Section \ref{sec:conc}.

\section{Definitions of matrix polynomials and notation}
\label{sec:Notation}

Recall that a matrix polynomial $P(\lambda) \in \RR[\lambda]^{n \times n}$ as in \eqref{equ:MatPol} is said to be \emph{$T$-even} if
\begin{equation} P(\lambda)^T = \sum_{k=0}^m \lambda^k P_k^T = P(-\lambda) \label{equ:MatrixPolynomial} \end{equation}
holds. In turn, \emph{$T$-odd} matrix polynomials are characterized by $P(\lambda)^T = - P(- \lambda)$. As mentioned in Section \ref{sec:intro}, it is easily seen that $P(\lambda)$ in \eqref{equ:MatrixPolynomial} is $T$-even if and only if $P_k^T = P_k$ holds for all matrix coefficients $P_k$ with even index $k \geq 0$ while $P_k^T = -P_k$ holds whenever the index $k$ is odd. The converse is true for $T$-odd matrix polynomials. Both structures have already been analyzed in \cite[Sec.\,6]{Mack_Diss}. The classes of regular, singular and unimodular matrix polynomials are defined as follows:

\begin{compactenum}
\item A matrix polynomial $P(\lambda) \in \RR[\lambda]^{n \times n}$ is called \emph{regular} if $\det(P(\lambda)) \neq 0$  and \emph{singular} otherwise (notice that $\det(P(\lambda)) \in \RR[\lambda]$).
\item A matrix polynomial $Q(\lambda) \in \RR[\lambda]^{n \times n}$ is called \emph{unimodular} if $\det(Q(\lambda))$ is a nonzero constant.
\end{compactenum}

Let $P(\lambda) \in \RR[\lambda]^{n \times n}$ be regular. We call $\mu \in \CC$ a (finite) \emph{eigenvalue} of $P(\lambda)$ if
$$P(\mu) = \sum_{k=0}^m \mu^k P_k \in \CC^{n \times n}$$
is a singular matrix.  Thus, $\mu \in \CC$ is an eigenvalue of $P(\lambda)$ if and only if $\det(P(\mu)) = 0$. Therefore, the set of all finite eigenvalues of $P(\lambda)$ coincides with the roots of $\det(P(\lambda)) \in \RR[\lambda]$ \cite[Sec.\,2]{MMMM06}. The \emph{algebraic multiplicity} of $\mu$ is defined as the multiplicity of $\mu$ as a root of $\det(P(\lambda))$. In addition, if $\mu \in \CC$ is some eigenvalue of $P(\lambda)$, the corresponding nullspace $\textnormal{null}(P(\mu))$ is called the \emph{eigenspace} for $\mu$. Its dimension is referred to as the \emph{geometric multiplicity} of $\mu$.

We define for any $d \geq \deg(P)$
$$ \textnormal{rev}_d P(\lambda) := \lambda^d P(\lambda^{-1}).$$
Then $\textnormal{rev}_d P(\lambda)$ is again a matrix polynomial, i.e. $\textnormal{rev}_d P(\lambda) \in \RR[\lambda]^{n \times n}$. It is called the \emph{$d$-reversal} corresponding to $P(\lambda)$ \cite[Def.\,2.2]{MMMM06}. In case $d = \textnormal{deg}(P)$, we call $\textnormal{rev}\,P(\lambda) := \textnormal{rev}_d\,P(\lambda)$ the reversal of $P(\lambda)$. It is easily verified that the finite eigenvalues of $\textnormal{rev}\,P(\lambda)$ are the reciprocals of the eigenvalues of $P(\lambda)$. In accordance with this observation, we call $\infty$ an eigenvalue of $P(\lambda)$ if zero in an eigenvalue of $\textnormal{rev}\,P(\lambda)$. The algebraic and geometric multiplicities of the eigenvalue $\infty$ are defined in terms of $\textnormal{rev}\,P(\lambda)$ and its finite eigenvalue $\mu = 0$  \cite[Def.\,2.3]{MMMM06}. The set of all eigenvalues of $P(\lambda) \in \RR[\lambda]^{n \times n}$ is called the \emph{spectrum} of $P(\lambda)$ and is denoted by $\sigma(P)$.

The following important property is intrinsic for the eigenvalues of $T$-even matrix polynomials $P(\lambda)$:

\begin{proposition} \label{prop:SpecStruct}
The spectrum $\sigma(P)$ of real, $T$-even matrix polynomials $P(\lambda)$ has a Hamiltonian structure. That is, $\sigma(P)$ is symmetric with respect to both the real and the imaginary axis.
\end{proposition}

Two matrix polynomials $S(\lambda), P(\lambda) \in \RR[\lambda]^{n \times n}$ are called unimodular \emph{equivalent}, if there exist unimodular matrix polynomials $U(\lambda), V(\lambda) \in \RR[\lambda]^{n \times n}$ such that $S(\lambda) = U(\lambda) P(\lambda) V(\lambda)$ holds. Linearizations for matrix polynomials are defined by unimodular equivalence as follows:

\begin{definition}[{Linearization, \cite[Def.\,2.12]{Dop18}}]  \label{def:linearization}
 Let $P(\lambda) \in \RR[\lambda]^{n \times n}$.
\begin{compactenum}
\item[$(i)$] Any matrix polynomial $\mathcal{L}(\lambda) = \lambda X + Y$ that can be expressed as \index{linearization}
\begin{equation}
\mathcal{L}(\lambda)  = U(\lambda) \begin{bmatrix} I_s & 0 \\ 0 &  P(\lambda) \end{bmatrix} V(\lambda) \in \RR[\lambda]^{(n+s) \times (n+s)} \label{equ:linearization}
\end{equation}
for two unimodular matrix polynomials $ U(\lambda), V(\lambda)$ of size $(n+s) \times (n+s)$ and  some $s \in \mathbb{N}_0$ is called a \textit{linearization} for $P(\lambda)$.
\item[$(ii)$] Assume $\textnormal{deg}(P) = k$. A linearization $\mathcal{L}(\lambda)$ for $P(\lambda)$ as in \eqref{equ:linearization} is called \textit{strong} (linearization) whenever $\textnormal{rev}_1 \mathcal{L}(\lambda)$ is a linearization for $\textnormal{rev}_k P(\lambda) = \textnormal{rev} P(\lambda),$ too.
\end{compactenum}
\end{definition}

Notice that unimodular matrix polynomials do not have any finite eigenvalues. Therefore, any linearization $\mathcal{L}(\lambda)$ as in \eqref{equ:linearization} of $P(\lambda)$ has the same finite eigenvalues (with the same algebraic and geometric multiplicities) as $P(\lambda)$ \cite{Dop18}. Furthermore, if $\mathcal{L}(\lambda)$ is strong, the same holds for the eigenvalue $\infty$ in case $\infty \in \sigma(P)$.

The problem of finding linearizations for matrix polynomials $P(\lambda) \in \RR[\lambda]^{n \times n}$ has been addressed in, e.g., \cite{MMMM06, PS_Diss}. Particular research has been done on conditioning \cite{HighMack06}, structure-preservation \cite{HighMack06-Sym} and nonstandard polynomial bases \cite{FaSa17, MackPer16, MackPer18}. In \cite{Dop18}, a new class of linearizations was introduced (so called block minimal bases linearizations) that has recently attracted much attention. The linearization we present in Theorem \ref{thm:linearization} will belong to this class.

Due to Proposition \ref{prop:SpecStruct}, we are particularly interested in $T$-even linearizations (whenever $P(\lambda) \in \RR[\lambda]^{n \times n}$ is $T$-even) to preserve the symmetries inherent to the spectrum of $P(\lambda)$. The structure of the $T$-even linearization $\mathcal{L}_P(\lambda)$ we define in \eqref{equ:Linearization} varies slightly depending on the parity of $\textnormal{deg}(P)$ (which can be even or odd).
Thus we define $M_P(\lambda)$ in Definition \ref{def:MP} depending on the degree of $P(\lambda)$ to treat both cases in Theorem \ref{thm:linearization} in a common framework.
Here and hereafter, we use the notation $\langle x,y \rangle$ to represent the scalar product $x^Ty \in \mathbb{R}$ of two vectors $x$ and $y$ and $\oplus$ to denote the direct sum of matrices, i.e. $A \oplus B = \textnormal{diag}(A,B)$ for any $A,B \in \mathbb{R}^{n \times n}$.

\begin{definition} \label{def:MP}
Assume $P(\lambda) \in \RR[\lambda]^{n \times n}$ is given as in \eqref{equ:MatrixPolynomial}.
\begin{enumerate}[(a)]
\item If $\textnormal{deg}(P)$ is odd, we define
 \begin{equation} M_P(\lambda) := \bigoplus_{k=0}^{\ell-1} (-1)^k \big(\lambda P_{d-2k} + P_{d-2k-1} \big) \in \mathbb{R}[\lambda]^{\ell n \times \ell n } \label{equ:bodyMP} \end{equation}
 with $d = \textnormal{deg}(P)$ and $\ell = (d+1)/2$.
\item If $\textnormal{deg}(P)$ is even, we define $M_P(\lambda)$ as in \eqref{equ:bodyMP} above with $d = \textnormal{deg}(P)+1$, $\ell = (d+1)/2$  and $P_d := 0_{n \times n}$.
\end{enumerate}
\end{definition}
Notice that $(\lambda P_{d-2k} + P_{d-2k-1})^T = - \lambda P_{d-2k} + P_{d-2k-1}$ holds for all summands in \eqref{equ:bodyMP} regardless of the parity of $\textnormal{deg}(P)$. That means $M_P(\lambda)^T = M_P(- \lambda)$, so $M_P(\lambda)$ is always $T$-even. With the definition $\Lambda_k(\lambda) := [ \, \lambda^k \; \lambda^{k-1} \; \cdots \; \lambda \; 1 \, ] \in \RR[\lambda]^{1 \times (k+1)}$ for any $k \geq 1$, we make the following important observation.

\begin{remark} \label{rem:body-equation}
According to the construction of $M_P(\lambda)$ for $P(\lambda)$ as in \eqref{equ:bodyMP} it can be verified by a direct calculation that
$$\big( \Lambda_\ell(-\lambda) \otimes I_n \big) M_P(\lambda) \big( \Lambda_{\ell}(\lambda)^T \otimes I_n \big) = P(\lambda)$$
holds. This property will be exploited in Section \ref{sec:LinearSystems}.
\end{remark}

Let $P(\lambda) \in \RR[\lambda]^{n \times n}$ be $T$-even. With the use of $M_P(\lambda)$ and
\begin{equation} L_k(\lambda) := \begin{bmatrix} 1 & -\lambda & & & \\ & 1 & -\lambda & & \\ & & \ddots & \ddots & \\ & & & 1 & -\lambda \end{bmatrix} \in \RR[\lambda]^{k \times (k+1)}, \quad k \geq 1, \label{equ:Lmatrix} \end{equation}
we present a structure-preserving, i.e. $T$-even, linearization $\mathcal{L}_P(\lambda)$ for $P(\lambda)$ in the following Theorem \ref{thm:linearization}. It is a block minimal bases pencil (as introduced in \cite{Dop18}) and was already used in \cite{FaSa18}. In particular, Theorem 3.3 in \cite{Dop18} applies to the matrix pencil $\mathcal{L}_P(\lambda)$ defined in \eqref{equ:Linearization} below and confirms that it is in fact a linearization for $P(\lambda)$.

\begin{theorem} \label{thm:linearization}
Let $P(\lambda) \in \RR[\lambda]^{n \times n}$ be $T$-even. Then the matrix pencil
\begin{equation}
\mathcal{L}_P(\lambda) := \left[ \begin{array}{c|c} M_P(\lambda) & L_{\ell-1}(-\lambda)^T \otimes I_n \\ \hline L_{\ell-1}(\lambda) \otimes I_n & 0 \end{array} \right] \in \mathbb{R}[\lambda]^{dn \times dn}
\label{equ:Linearization}
\end{equation}
defined for $P(\lambda)$ with $M_P(\lambda)$, $d$ and $\ell$ given as in Definition \ref{def:MP} and $L_{\ell-1}(\lambda)$ as introduced in \eqref{equ:Lmatrix} is $T$-even. Moreover, $\mathcal{L}_P(\lambda)$ is a strong linearization for $P(\lambda)$ if $\textnormal{deg}(P)$ is odd and a linearization for $P(\lambda)$ if $\textnormal{deg}(P)$ is even.
\end{theorem}

Due to the linearization property, the matrix pencil $\mathcal{L}_P(\lambda) \in \RR[\lambda]^{dn \times dn}$ defined in \eqref{equ:Linearization} has the same finite eigenvalues as $P(\lambda) \in \RR[\lambda]^{n \times n}$ (from whose matrix coefficients it is defined). Moreover, since $\mathcal{L}_P(\lambda)$ is $T$-even whenever $P(\lambda)$ is $T$-even, we call $\mathcal{L}_P(\lambda)$ a structure-preserving linearization for $P(\lambda)$.
To illustrate the form of $\mathcal{L}_P(\lambda)$ consider the following example.

\begin{example}
The linearization $\mathcal{L}_P(\lambda)$ defined for a $T$-even matrix polynomial $P(\lambda) \in \RR[\lambda]^{n \times n}$ has a very sparse and clear structure. This is illustrated below for $P(\lambda) = \sum_{k=0}^7 \lambda^k P_k$ of degree seven. Writing $\mathcal{L}_P(\lambda)$ in the form $\lambda X + Y$ for two $7n \times 7n$ matrices $X$ and $Y$ we have
\begin{equation} \begin{aligned} \mathcal{L}_P(\lambda) &= \left[ \begin{array}{cccc|ccc} -P_7 & & & & & & \\ & P_5 & & & I_n & & \\ & &-P_3 & & & I_n & \\  & & & P_1 & & &I_n \\ \hline &-I_n & & & & & \\ & & -I_n & & & & \\ & & & -I_n & & & \end{array} \right]
\lambda \\ & \hspace{3cm} + \left[ \begin{array}{cccc|ccc} -P_6 & & & & I_n & & \\ & P_4 & & & & I_n & \\ & & -P_2 & & & & I_n \\ & & & P_0 & & & \\ \hline  I_n & & & & & & \\ & I_n & & & & & \\ & & I_n & & & & \end{array} \right].
\end{aligned} \label{equ:example-Lin} \end{equation}
Since $P(\lambda)$ was assumed to be $T$-even, it is seen directly that $Y$ is symmetric while $X$ is skew-symmetric.
In addition, notice that, if $P(\lambda)$ was only of degree six, $\mathcal{L}_P(\lambda)$ as defined in \eqref{equ:Linearization} would be as in \eqref{equ:example-Lin} with $P_7=0$.
\end{example}

In general, determining $\sigma(P)$ for a matrix polynomial $P(\lambda)$ is sometimes referred to as the \emph{polynomial eigenvalue problem} (PEP). If $P(\lambda) = \lambda X + Y$ is a matrix pencil, the term generalized eigenvalue problem (GEP) is often used. A common way to solve a PEP corresponding to $P(\lambda) \in \RR[\lambda]^{n \times n}$ is to compute $\sigma(P)$ (or just a part of it) through a linearization $\mathcal{L}(\lambda)$ for $P(\lambda)$ using a method for GEPs. Notice that the size of $\mathcal{L}(\lambda)$ is usually much larger than the size of $P(\lambda)$ (depending on the degree of $P(\lambda)$). Thus it is often appropriate (or even necessary) not to compute all eigenvalues of $\mathcal{L}(\lambda)$ but only some (e.g. in a predefined area of the complex plane).  For such purposes Krylov subspace methods are among the most appropriate algorithms (cf. \cite{Templates} for an overview of different Krylov subspace algorithms). Hereby, the area where eigenvalues are to be found is controlled via a shift parameter $\zeta \in \CC$. With some abuse of terminology, a Krylov subspace method can be called rational if it admits changes of this shift parameter during its iteration, see \cite{Ruhe84,Ruhe98}.

All subsequent investigations mainly aim for the construction of a rational Krylov subspace algorithm to determine eigenvalues of $\mathcal{L}_P(\lambda)$ defined as in \eqref{equ:Linearization} for some given $T$-even matrix polynomial $P(\lambda) \in \RR[\lambda]^{n \times n}$. Hereby, the $T$-even structure of $\mathcal{L}_P(\lambda)$ is exploited to preserve the spectral symmetries.

\section{The Even-IRA algorithm}
\label{sec:EvenIRA}

According to Proposition \ref{prop:SpecStruct}, the spectrum of a (real) $T$-even matrix polynomial is symmetric with respect to the real and imaginary axis.
Numerical algorithms respecting this spectral symmetry will in general be more accurate than standard methods \cite{MackVibr}. In addition, numerical methods that ignore the special structure may often produce (physically) less meaningful results \cite{MeerTiss}. Therefore, our focus in the development of a reliable eigensolver for $T$-even polynomial eigenvalue problems is twofold: on the one hand on the application of a structure-preserving linearization (see Theorem \ref{thm:linearization}) and on the other hand on a method that profitably exploits this structure. One method taking the $T$-even structure into account is the Even-IRA algorithm presented in \cite{EvenIRA}.
It belongs to the class of Krylov subspace methods and is a sophisticated variant of the Krylov-Schur algorithm \cite{KrylovSchur} customized for real $T$-even generalized eigenvalue problems.
Other methods for solving $T$-even polynomial eigenvalue problems can be found in, e.g., \cite{Bassour20, MehrWat01}.

The \textsc{Even-IRA} algorithm is designed to determine a part of the spectrum of a regular $T$-even matrix pencil $\mathcal{G}(\lambda) = \lambda X + Y \in \RR[\lambda]^{m \times m}$ close to a predefined target in the complex plane. To preserve the Hamiltonian eigenvalue structure, a special spectral transformation is applied to preserve the $\pm$ matching pairs of eigenvalues.

In particular, whenever $\mathcal{G}(\lambda) = \lambda X + Y$ is regular and $T$-even, i.e. $X=-X^T$ and $Y=Y^T$ holds, and some shift $\zeta \notin \sigma(\mathcal{G})$ is given in a region of the complex plane where eigenvalues are to be found, then in \cite{EvenIRA} the transformation
\begin{equation}  \mathcal{G}(\zeta) = \zeta X + Y \mapsto K(\zeta) = \mathcal{G}(\zeta)^{-T}X \mathcal{G}(\zeta)^{-1}X \in \CC^{m \times m} \label{equ:SpecTransform} \end{equation}
is considered. Notice that a similar spectral transformation already appeared in \cite{BenFaStoll,MehrWat01,Wat04} in the context of skew-Hamiltonian/Hamiltonian eigenvalue problems and the symplectic Lanczos process.  Whenever $\mathcal{G}(\mu)x=0$ holds for some $\mu \in \CC$ and $x \in \CC^m$, it is easily confirmed that $K(\zeta)x = \theta x$ follows, where $\theta = (\mu^2 - \zeta^2)^{-1}$. Thus, any two finite eigenvalues $\mu$ and $-\mu$ of $\mathcal{G}(\lambda)$ are mapped to the same eigenvalue $\theta \in \sigma(K(\zeta))$. Due to this fact, $\pm$ matching pairs of eigenvalues are preserved. On the other hand, all eigenvalues of $K(\zeta)$ are necessarily of even multiplicity. Notice the following two important facts: 
\begin{itemize}
\item Whenever some eigenvalue $\theta \in \sigma(K(\zeta))$ has been found, it gives rise to a $\pm$ matching
pair of two eigenvalues of $\mathcal{G}(\lambda)$, namely
\begin{equation} \mu = \sqrt{(1/\theta) + \zeta^2} \quad \textnormal{and} \quad \widehat{\mu} = - \sqrt{(1/\theta) + \zeta^2}. \label{equ:rev-transform} \end{equation}
\item The matrix $K(\zeta)$ from \eqref{equ:SpecTransform} will in general be complex but remains real whenever $\zeta \in \RR$ or $\zeta \in i \RR$. In case $\zeta = a + bi$ with nonzero real and imaginary parts, a slightly different spectral transformation can be considered, see \cite[Rem.\,2.1]{EvenIRA}, to stay in real arithmetic. 
\end{itemize}

In \cite{EvenIRA} the authors suggest to apply the implicitly restarted Krylov-Schur method \cite{KrylovSchur} to the matrix $K(\zeta)$ in \eqref{equ:SpecTransform} to find some, say $s \in \mathbb{N}$,  eigenvalues of $\mathcal{G}(\lambda)$. That is, if $v_1, \ldots , v_s$ is an orthonormal basis of the Krylov space
\begin{equation} \mathcal{K}_s(K(\zeta),x) = \textnormal{span} \lbrace x, K(\zeta)x, K(\zeta)^2x, \ldots , K(\zeta)^{s-1}x \rbrace \label{equ:Krylovspace} \end{equation} for some $x \in \RR^m$ (computed by the Arnoldi method, see \cite{Templates})
and $V = [ \, v_1 \; \cdots \; v_s \, ] \in \RR^{m \times s}$, in general some of the eigenvalues of $K(\zeta)$ of largest magnitude are well approximated by some of the $s$ eigenvalues of $V^TK(\zeta)V$.
This process can now be (implicitly) restarted using the Krylov-Schur restart strategy \cite[Sec.\,3]{KrylovSchur} until all $s$ eigenvalues of $V^TK(\zeta)V$ serve as good approximations to eigenvalues of $K(\zeta)$. This approach is called the \textsc{Even-IRA} algorithm (details on the practical implementation of the algorithm can be found in \cite[Sec.\,4]{EvenIRA}). Additional information on how eigenvectors may be captured can be found in \cite[p.\,4074ff]{EvenIRA}.

The basis of our algorithm is the \textsc{Even-IRA} algorithm. As this method is designed for $T$-even matrix pencils, it cannot be used directly for $T$-even matrix polyomials $P(\lambda)$ of degree $>1$. To solve the PEP for $P(\lambda)$, we apply the \textsc{Even-IRA} algorithm to the structure-preserving linearization $\mathcal{L}_P(\lambda)$ from \eqref{equ:Linearization}. The sparse block structure of $\mathcal{L}_P(\lambda)$ turns out to be very beneficial for the computation of matrix-vector-products $K(\zeta)x$ (which are necessary to build the Krylov space, see \eqref{equ:Krylovspace}). In fact, we show in Section \ref{sec:LinearSystems} that $K(\zeta)x$ can be computed in a cheap and reliable way without ever forming $K(\zeta)$ and $\mathcal{L}_P(\zeta)$ explicitly. In Section \ref{sec:RatArnoldiDecomp}, we will modify the \textsc{Even-IRA} algorithm so that it is able to handle changes of the shift parameter $\zeta$ during the Arnoldi iteration and the restart process. This makes it possible to accelerate convergence or to control/change the regions in the complex plane where eigenvalues are to be found.

\section{The efficient computation of matrix-vector-products $K(\zeta)x$}
\label{sec:LinearSystems}

Assume that $P(\lambda) \in \RR[\lambda]^{n \times n}$ is some $T$-even matrix polynomial and let $\mathcal{L}_P(\lambda) \in \RR[\lambda]^{dn \times dn}$ be defined as in \eqref{equ:Linearization}. Recall that $P(\lambda)$ and $\mathcal{L}_P(\lambda)$ share the same finite eigenvalues. As outlined in Theorem \ref{thm:linearization}, $\mathcal{L}_P(\lambda)$ is $T$-even, so the \textsc{Even-IRA} algorithm can be applied to $\mathcal{L}_P(\lambda)$ to determine a part of the finite spectrum of $P(\lambda)$. In consideration of large-scale-problems\footnote{Recall that the size of $\mathcal{L}_P(\lambda)$ is (depending on the degree of $P(\lambda)$) much larger than the size of $P(\lambda)$. A major part of the computational cost that is raised by a Krylov subspace method such as the \textsc{Even-IRA} algorithm usually comes from the computation of the matrix-vector-products to form the Krylov space (see \eqref{equ:Krylovspace}).
Thus, to achieve a reasonable efficiency of our method, it is necessary to guarantee the fast and cheap computation of matrix-vector-products $K(\zeta)v$, where $K(\zeta) \in \CC^{dn \times dn}$ is the matrix defined for $\mathcal{L}_P(\lambda)$ in \eqref{equ:SpecTransform}, $\zeta \in \CC$ is some shift not contained in the spectrum of $P(\lambda)$, and $v \in \CC^{dn}$.}, the sparsity and structure of $\mathcal{L}_P(\lambda)$ can be exploited to significantly increase the computational speed in calculating $K(\zeta)v$. This effective computational approach in explained in detail in this section (see also  \cite{FaSa18} and \cite[Sec.\,5.2]{PS_Diss}).

\begin{remark}
As we are only considering polynomial eigenvalue problems given by real $T$-even matrix polynomials, matrix-vector-multiplications $K(\zeta)v$ will only involve real vectors $v \in \RR^{dn}$ in all subsequent sections (even if $\zeta$ and $K(\zeta)$ are complex). However, the technique to perform matrix-vector-multiplications is valid even if $v \in \CC^{dn}$, so we give a general treatment here.
\end{remark}

To begin, assume $\zeta \in \CC$ is not contained in $\sigma(P)$ and let $v \in \CC^{dn}$ be given. Moreover, let $\mathcal{L}_P(\lambda) = \lambda X + Y$ as in \eqref{equ:Linearization}. Explicitly, the matrix-vector-product $K(\zeta)v$ can be written as
\begin{equation} K(\zeta)v = \left( \mathcal{L}_P(\zeta)^{-T}X \mathcal{L}_P(\zeta)^{-1}X \right) v. \label{equ:MatrVecProd} \end{equation}
Actually, \eqref{equ:MatrVecProd} can be evaluated using four consecutive matrix-vector-mul\-ti\-pli\-cations. The matrix-vector-products with $X$, where $X \in \RR^{dn \times dn}$, can be evaluated directly and quickly by exploiting the sparsity of $X$. Moreover, as $X$ has a clear and determined block-structure, a matrix-vector-multiplication $Xv$ can entirely be carried out implicitly, that is, without forming $X$ at all, on its nonzero $n \times n$ blocks. The matrix-vector-products with $\mathcal{L}_P(\zeta)^{-1}$ and $\mathcal{L}_P(\zeta)^{-T} = \mathcal{L}_P(-\lambda)^{-1}$ can be realized by solving linear systems with $\mathcal{L}_P(\zeta)$ and $\mathcal{L}_P(- \zeta)$, respectively. However, the size of both matrices is $dn \times dn$ and, therefore, can be rather large. Fortunately, a linear-systems-solve with $\mathcal{L}_P(\zeta)$ can be traced back to solely $n \times n$ computations. The solution of a linear system with $\mathcal{L}_P(\lambda)$ can essentially be reduced to the solution of a linear system involving $P(\zeta) \in \CC^{n \times n}$. This provides an economic approach for the determination of these products since, for instance, the computational cost of an LU decomposition for $\mathcal{L}_P(\zeta)$ is within $\mathcal{O}(d^3n^3)$ while it is only $\mathcal{O}(n^3)$ for $P(\zeta)$ if no sparsity patterns are taken into account. For sparse matrices the cost is about $\mathcal{O}(d \cdot \textnormal{nz} )$ and $\mathcal{O}( \textnormal{nz})$, respectively, where nz denotes the number of nonzero entries. Moreover, the storage requirements for the LU factors for $P(\zeta)$ are way below those for the LU factors of $\mathcal{L}_P(\lambda)$. 

Assume that $\mathcal{L}_P(\zeta)^{-1}v$ is to be computed, i.e. the linear system
\begin{equation}  \mathcal{L}_P(\zeta)y = \left[ \begin{array}{c|c} M_P(\zeta) & L_{\ell-1}(-\zeta)^T \otimes I_n \\ \hline L_{\ell-1}(\zeta) \otimes I_n & 0 \end{array} \right] \begin{bmatrix} y_1 \\ y_2 \end{bmatrix} = \begin{bmatrix} x_1 \\ x_2 \end{bmatrix} = x \label{equ:System0} \end{equation}
is to be solved for a given vector $x \in \CC^{dn}$. Let $y^\star$ be the solution of \eqref{equ:System0} which is unique since $\mathcal{L}_P(\zeta)$ is nonsingular (due to the fact that $\zeta \notin \sigma(P)$).
Let $x,y \in \CC^{dn}$ be partitioned as $x_1, y_1 \in \CC^{\ell n}$ and $x_2, y_2 \in \CC^{(\ell -1)n}$ and assume that $y^\star = [ \, (y_1^\star)^T \; (y_2^\star)^T \, ]^T$ is partitioned accordingly.
The structure of $\mathcal{L}_P(\zeta)$ reveals that \eqref{equ:System0} can be rewritten as a system of two equations for the unknown vectors $y_1$ and $y_2$ as
\begin{align}
M_P(\zeta)y_1 + \left( L_{\ell-1}(-\zeta)^T \otimes I_n \right)y_2 &= x_1 \quad \textnormal{and} \label{equ:System1}\\
\left( L_{\ell-1}(\zeta) \otimes I_n \right)y_1 &= x_2. \label{equ:System2}
\end{align}
Notice that \eqref{equ:System2} is an underdetermined system with $L_{\ell-1}(\zeta) \otimes I_n \in \CC^{(\ell -1)n \times \ell n}$. Moreover, $\rank(L_{\ell-1}(\zeta) \otimes I_n)) = (\ell-1)n$ holds regardless of the choice of $\zeta$. Therefore, the nullspace of $L_{\ell-1}(\zeta) \otimes I_n$ is always $n$-dimensional and easily determined since $$ \big( L_{\ell}(\zeta) \otimes I_n \big) \big( \Lambda_{\ell}(\zeta)^T \otimes I_n \big) = 0_{(\ell-1)n \times n}.$$
Therefore we have $\nullspace(L_{\ell}(\zeta) \otimes I_n) = \lbrace (\Lambda_{\ell}(\zeta)^T \otimes I_n)r \,; \, r \in \CC^n \rbrace$.
Consequently, any solution $y_1$ for \eqref{equ:System2} has the form $y_1 = \widehat{y}_1 + (\Lambda_\ell(\zeta)^T \otimes I_n)r$ where $\widehat{y}_1 \in \CC^{\ell n}$ solves \eqref{equ:System2} and $r \in \CC^n$ is arbitrary (i.e. $(\Lambda_\ell(\zeta)^T \otimes I_n)r$ is a solution to the homogenous system corresponding to \eqref{equ:System2}). In fact, once some particular solution $\widehat{y}_1$ has been found, there exists some unique $r^\star$ such that $y_1^\star = \widehat{y}_1 + (\Lambda_{\ell}(\zeta)^T \otimes I_n)r^\star \in \CC^{\ell n}$. With this characterization of $y_1^\star$ at hand, it follows from \eqref{equ:System1} that
\begin{equation}
M_P(\zeta) \big( \widehat{y}_1 + (\Lambda_{\ell}(\zeta)^T \otimes I_n)r^\star \big) + \left( L_{\ell}(-\zeta)^T \otimes I_n \right)y_2^\star = x_1
\label{equ:System3}
\end{equation}
holds. Multiplying \eqref{equ:System3} by $\Lambda_\ell(-\zeta) \otimes I_n$ from the left eliminates the second term since $(\Lambda_{\ell}(-\zeta) \otimes I_n) (L_\ell(-\zeta)^T \otimes I_n) = 0$. After some reordering we obtain from \eqref{equ:System3}
\begin{equation}
\underbrace{\big( \Lambda_\ell(-\zeta) \otimes I_n \big) M_P(\zeta) \big( \Lambda_{\ell}(\zeta)^T \otimes I_n \big)}_{= \, P(\zeta)}r^\star = x_1 - \big( \Lambda_\ell(-\zeta) \otimes I_n \big) M_P(\zeta) \widehat{y}_1.
\label{equ:System4}
\end{equation}
Notice that the left-hand-side of \eqref{equ:System4} simplifies to $P(\zeta)r^\star$ in accordance with Remark \ref{rem:body-equation}. In addition,
as $P(\zeta)$ is nonsingular, $r^\star$ is the unique solution of \eqref{equ:System4}. Thus, in other words, for any fixed particular solution $\widehat{y}_1$ of \eqref{equ:System2}, the unique solution $r^\star$ of the $n \times n$ linear system
\begin{equation}
P(\zeta)r = x_1 - \big( \Lambda_\ell(-\zeta) \otimes I_n \big) M_P(\zeta) \widehat{y}_1
\label{equ:System5}
\end{equation}
determines the first part  $y_1^\star  =  \widehat{y}_1 + (\Lambda_{\ell}(\zeta)^T \otimes I_n)r^\star \in \CC^{\ell n}$ of the solution vector $y^\star$. Once $y_1^\star$ has been found, $y_2^\star \in \CC^{(\ell -1)n}$ will be the unique solution of the overdetermined system
\begin{equation}
\left( L_{\ell-1}(-\zeta) \otimes I_n \right)y_2 = x_1 - M_P(\zeta) y_1^\star
\label{equ:System6}
\end{equation}
since \eqref{equ:System1} and \eqref{equ:System2} are satisfied if and only if $y_1 = y_1^\star$ and $y_2 = y_2^\star$. The computations of a particular solution $\widehat{y}_1$ of \eqref{equ:System2} and the solution $y_2^\star$ of \eqref{equ:System6} can be carried out by forward and backward substitution and both require $\mathcal{O}(\ell n)$ flops. In particular:

\begin{compactenum}

\item A solution $\widehat{y}_1 \in \CC^{\ell n}$ for \eqref{equ:System2}, i.e.
\begin{equation} \begin{bmatrix} I_n & - \zeta I_n & & & \\ & I_n & -\zeta I_n & & \\ & & \ddots & \ddots & \\ & & & I_n & -\zeta I_n \end{bmatrix} \begin{bmatrix} y_{1,1} \\ y_{1,2} \\ \vdots \\ y_{1,\ell-1} \\ y_{1, \ell} \end{bmatrix} = \begin{bmatrix} v_{1,1} \\ v_{1,2} \\ \vdots \\ v_{1,\ell-1} \end{bmatrix}, \quad y_{1,k}, v_{1,k} \in \CC^n,
\label{equ:System7}
\end{equation}
can be found by backward substitution. If $\widehat{y}_1$ and $v_1 \in \CC^{(\ell -1)n}$ are partitioned as in \eqref{equ:System7} and $\widehat{y}_{1,\ell} =0$ is chosen, then $\widehat{y}_{1,1}, \ldots , \widehat{y}_{1, \ell -1}$ are uniquely determined through the recurrence relation
$\widehat{y}_{1,k} = v_{1,k} + \zeta \widehat{y}_{1,k+1}$ for $k=\ell-1, \ldots , 1.$

\item The unique solution $y_2^\star$ of \eqref{equ:System6}, i.e.
\begin{equation} \begin{bmatrix} I_n & & & \\ \zeta I_n & I_n & & \\ & \zeta I_n & \ddots & \\ & & \ddots & I_n \\ & & & \zeta I_n \end{bmatrix} \begin{bmatrix} y_{2,1} \\ y_{2,2} \\ \vdots \\ y_{2, \ell-1} \end{bmatrix} = x_1 - M_P(\zeta) y_1^\star =: \begin{bmatrix} w_1 \\ w_2 \\ \vdots \\ w_{\ell-1} \\ w_\ell \end{bmatrix}, \; y_{2,k}, w_k \in \CC^n,
\label{equ:System8}
\end{equation}
can be found by forward substitution. If $y_2^\star$ is partitioned as $y_2$ in \eqref{equ:System8}, then $y_{2,1}^\star = w_1$ and $y_{2,2}^\star, \ldots , y_{2,\ell}^\star$ are uniquely determined by the recurrence
$ y_{2,k}^\star = w_k - \zeta y_{2,k-1}^\star$ for $k=2, \ldots , \ell.$
\end{compactenum}

\begin{remark} \label{rem:overdetermined}
Notice that, although \eqref{equ:System8} is an overdetermined system for $y_2 \in \CC^{(\ell -1)n}$ which usually need not have a solution, there is a unique solution $y_2^\star$ for \eqref{equ:System8} since we assumed $\mathcal{L}_P(\zeta)y = x$ to be uniquely solvable.
\end{remark}
For the determination of matrix-vector-products $\mathcal{L}_P(\zeta)^{-T}v$, the $T$-even structure of $\mathcal{L}_P(\zeta)$ can be exploited. In particular, $\mathcal{L}_P(\zeta)^{-T} = (\mathcal{L}_P(\zeta)^T)^{-1} = \mathcal{L}_P(-\zeta)^{-1}$. In order to find $\mathcal{L}_P(-\zeta)^{-1}v$, the same approach as above can be used involving $- \zeta$ instead of $\zeta$. In particular,  in \eqref{equ:System5} the matrix $P(-\zeta)$ instead of $P(\zeta)$ will show up. If an LU decomposition $P(\zeta)=LU$ has been computed to solve the linear system with $P(\zeta)$ in \eqref{equ:System5}, this factorization can be reused to solve the system with $P(- \zeta)$ since $P(- \zeta) = P(\zeta)^T = U^TL^T$.

In conclusion, the procedure described in this section presents an efficient way to calculate matrix-vector-products of the form \eqref{equ:MatrVecProd}. Whenever $d \ll n$, the complexity of the overall method is dominated by the cost of the LU decomposition of $P(\zeta)$ which is $\mathcal{O}(\textnormal{nz})$ or $\mathcal{O}(n^3)$ depending on whether $P(\zeta)$ is sparse or not.

\section{The rational Arnoldi decomposition}
\label{sec:RatArnoldiDecomp}

Let $P(\lambda) \in \RR[\lambda]^{n \times n}$ be some $T$-even matrix polynomial and let $\mathcal{L}_P(\lambda) = \lambda X + Y \in \RR[\lambda]^{dn \times dn}$ and
\begin{equation} K(\zeta) = \mathcal{L}_P(\zeta)^{-T}X \mathcal{L}_P(\zeta)^{-1}X = \mathcal{L}_P(-\zeta)^{-1}X \mathcal{L}_P(\zeta)^{-1}X, \quad \zeta \notin \sigma(P), \label{equ:K_Matrix} \end{equation}
be defined for $P(\lambda)$ as in \eqref{equ:Linearization} and \eqref{equ:SpecTransform}, respectively.

Recall from Section \ref{sec:EvenIRA} that $K(\zeta)$ as in \eqref{equ:K_Matrix} stays real whenever $\zeta$ is real or purely imaginary. We will assume for the moment that either of them holds to stay within real arithmetics. So, let $v_1 \in \RR^{dn}$ be some normalized vector and suppose that (for instance as part of the \textsc{Even-IRA} algorithm) $m \in \mathbb{N}$ steps of the Arnoldi process (cf. \cite[Alg.\,7.3]{Templates}) have been performed for $K(\zeta)$. That is, we are with an Arnoldi decomposition for $K(\zeta) \in \RR^{dn \times dn}$ of the form
\begin{equation} K(\zeta) \left[ \, \drawmatrix[width=.8]{V_m} \,  \right] = \left[ \, \drawmatrix[width=.8]{V_m} \,  \right] \left[ \, \drawmatrix[lower banded]{\;\,T_m} \, \right] + t_{m+1,m}v_{m+1}e_m^T
= \left[ \, \drawmatrix[width=.92]{V_{m+1}} \,  \right] \left[ \, \drawmatrix[size=1.1,width=1, lower banded]{\;\;\underline{T}_m} \, \right] \label{equ:ArnoldiDec0} \end{equation}
where $V_{m+1} = [ \, v_1 \; \cdots \; v_{m+1} \, ] = [ \; V_m \; v_{m+1} \, ] \in \RR^{dn \times (m+1)}$.
The following statements hold for the vectors and matrices involved in \eqref{equ:ArnoldiDec0}: \vspace{0.1cm}
\begin{compactenum}
\item The columns $v_1, \ldots , v_{m+1} \in \RR^{dn}$ of $V_{m+1}$ form an orthonormal basis of the Krylov space $\mathcal{K}_{m+1}(K(\zeta), v_1)$ where $v_1$ is the starting vector of the Arnoldi iteration.
\item The matrices $T_m =[t_{i,j}]_{i,j} \in \RR^{m \times m}$ and $\underline{T}_m := I_{m+1,m}T_m + t_{m+1,m}v_{m+1}e_m^T \in \RR^{(m+1) \times m}$ have upper-Hessenberg structure where $e_m$ denotes the $m$-th unit vector in $\RR^{m}$.
\end{compactenum} \vspace{0.1cm}
The eigenvalues of $T_{m} \in \RR^{m \times m}$ are called Ritz values of $K(\zeta)$ with respect to $\mathcal{K}_{m+1}(K(\zeta), v_1)$. According to the Rayleigh-Ritz principle (cf. \cite[Sec.\,7]{Demm97}), these values are used as approximations to the eigenvalues of $K(\zeta)$ by the \textsc{Even-IRA} algorithm (cf. \cite[pp.\,4074ff]{EvenIRA}, see also Section \ref{sec:EvenIRA}). Recall that the spectral transformation $\mu \mapsto (\mu^2 - \zeta^2)^{-1}$ corresponding to the transformation \eqref{equ:SpecTransform} causes eigenvalues $\mu$ of $\mathcal{L}_P(\lambda)$ close to $\zeta$ to be of large magnitude. In the first place, these will be well approximated by eigenvalues of $T_m$. Starting with the decomposition \eqref{equ:ArnoldiDec0}, the \textsc{Even-IRA} algorithm performs several Krylov-Schur restarts (see \cite{KrylovSchur}) until convergence of the desired number of eigenvalues was observed. As soon as $t_{m+1,m}$ in \eqref{equ:ArnoldiDec0} becomes zero, $K(\zeta)V_m = V_m T_m$ holds and all eigenvalues of $T_m$ are exact eigenvalues of $K(\zeta)$. Finally, the reverse transformation \eqref{equ:rev-transform} reveals eigenvalues of $\mathcal{L}_P(\lambda)$ close to $\zeta$.

Now notice that $K(\zeta)$ in \eqref{equ:K_Matrix} is nonsingular if and only if $X$ is nonsingular. Therefore, assuming $X$ to be nonsingular, we have
\begin{equation} \begin{aligned}
K(\zeta)^{-1} &= X^{-1} \mathcal{L}_P(\zeta) X^{-1} \mathcal{L}_P(- \zeta) \vphantom{\big( \big)}
              = X^{-1} \big( \zeta X + Y \big) X^{-1} \big( - \zeta X + Y \big) \\ &= - \zeta^2 I_{dn} + X^{-1}YX^{-1}Y
              = \big( X^{-1}Y \big)^2 - \zeta^2 I_{dn}
\end{aligned} \label{equ:KInverse} \end{equation}
so that $K(\zeta) = ( (X^{-1}Y)^2 - \zeta^2 I_{dn})^{-1}$. For all further considerations we let $G := (X^{-1}Y)$ whenever it exists so that $K(\zeta) = (G^2- \zeta^2 I_{dn})^{-1}$.

Now, whenever $X \in \RR^{dn \times dn}$ is nonsingular and $G$ exists, \eqref{equ:KInverse} can be taken into account and \eqref{equ:ArnoldiDec0} may be rewritten in terms of $G^2$ as
\begin{equation} \left[ \, \drawmatrix{G^2} \,  \right] \left[ \, \drawmatrix[width=.92]{V_{m+1}} \,  \right] \left[ \, \drawmatrix[size=1.2,width=1.1, lower banded]{\;\;\underline{T}_{m}} \, \right] = \left[ \, \drawmatrix[width=.92]{V_{m+1}} \,  \right] \left[ \, \drawmatrix[size=1.2,width=1.1, lower banded]{\;\;\underline{H}_{m}} \, \right]
\label{equ:ArnoldiDec1} \end{equation}
where $\underline{H}_{m} := \zeta^2 \underline{T}_{m} + I_{m+1,m} \in \RR^{(m+1) \times m}$ is again of upper-Hessenberg form. A decomposition of the form \eqref{equ:ArnoldiDec1} is called a generalized rational Arnoldi decomposition for $G^2$ in \cite[(1.1)]{BerlGuett15}, so we adapt this terminology here. When working with $K(\zeta) = (G^2- \zeta^2 I_{dn})^{-1}$ we will mainly consider rational decompositions as in \eqref{equ:ArnoldiDec1} instead of Arnoldi decompositions as in \eqref{equ:ArnoldiDec0} from now on.

Now regarding \eqref{equ:ArnoldiDec1}, the eigenvalues of the matrix pencil $\lambda T_{m} + H_{m} \in \RR[\lambda]^{m \times m}$ (where $T_m$ and $H_m$ are given by the first $m$ rows of $\underline{T}_m$ and $\underline{H}_m$, respectively) will in general be good approximations of the eigenvalues of $G^2$, see \cite{Ruhe98}. Certainly it holds that $\sigma(G^2) = \sigma(G)^2 = \sigma(\mathcal{L}_P)^2$, so the square roots $+ \sqrt{\theta}$ and $- \sqrt{\theta}$ of eigenvalues $\theta \in \CC$ found from $\lambda T_m + H_m$ approximate eigenvalues of $\mathcal{L}_P(\lambda)$ and, in turn, $P(\lambda)$. It is a crucial observation regarding the rational \textsc{Even-IRA} algorithm presented in Section \ref{sec:RatEvenIRA} that this relationship holds even if $G$ (and hence $G^2$) does not exist.

Now it is important to notice that, in \eqref{equ:ArnoldiDec1}, only $\underline{T}_m$ and $\underline{H}_m$ directly depend on $\zeta$ but $G^2$ does not. In comparison to the Arnoldi decomposition \eqref{equ:ArnoldiDec0} - where $K(\zeta)$ appears on the left-hand-side and explicitly depends on $\zeta$ - this enables us to extend the decomposition \eqref{equ:ArnoldiDec1} while changing the shift $\zeta$ to some newly chosen value $\xi \in \CC$. This is not possible for the standard Arnoldi decomposition \eqref{equ:ArnoldiDec0} and cannot be realized in the \textsc{Even-IRA} algorithm. In particular, instead of calculating and orthogonalizing $K(\zeta)v_{m+1}$ to extend \eqref{equ:ArnoldiDec1} (as in the Arnoldi iteration), we may use the vector $K(\xi)v_{m+1}$ for some new shift $\xi$.

\begin{remark}
We will assume throughout this section that $G \in \mathbb{R}^{dn \times dn}$ exists since this will be helpful to illustrate the forthcoming computations. This assumption will be dropped in the next section since the regularity of $X$ is actually not required to perform the rational \textsc{Even-IRA} algorithm outlined in Section \ref{sec:RatEvenIRA}.
\end{remark}

In Section \ref{ssec:extension} we show how the rational Arnoldi decomposition \eqref{equ:ArnoldiDec1} can be extended to increase the dimension of the underlying Krylov space (spanned by the columns of $V_{m+1}$).
To this end, we distinguish between the cases where either $\xi \in \RR$ or $\xi \in i \RR$ holds (only in these cases $\xi^2$ and $K(\xi)$ are real) and where $\xi = a+bi$ is complex with nonzero real and imaginary parts (which implies $K(\xi)$ to be non-real). In the second case, we may still remain in real arithmetics if the real and imaginary parts of $K(\xi)v_{m+1}$ are considered separately\footnote{Notice that the authors from \cite{EvenIRA} deal with complex shifts in another way by changing the definiton of $K(\zeta)$, see \cite[Rem.\,2.1]{EvenIRA}.}.

\subsection{The extension of a rational Arnoldi decomposition}
\label{ssec:extension}

Assume we are given a decomposition as in \eqref{equ:ArnoldiDec1} obtained from $m$ steps of the Arnoldi iteration applied to $K(\zeta)$ for some shift $\zeta \notin \sigma(\mathcal{L}_P)$. Now let $\xi \notin \sigma(\mathcal{L}_P)$ be some new shift parameter.
First assume that either $\xi \in \RR$ or $\xi \in i \RR$ holds, so $K(\xi)$ is a real matrix and $K(\xi)v_{m+1}$ is a real vector of size $dn$. 
The Gram-Schmidt-orthogonalization of $K(\xi)v_{m+1}$ against $v_1, \ldots , v_{m+1}$ yields
\begin{equation}
\widetilde{v}_{m+2} = K(\xi)v_{m+1} - \begin{bmatrix} v_1 & \cdots & v_{m+1} \end{bmatrix} \begin{bmatrix} t_{1,m+1} \\ \vdots \\ t_{m+1,m+1} \end{bmatrix},
 \label{equ:ArnoldiTransf}
\end{equation}
where $t_{i,m+1} = \langle K(\xi)v_{m+1}, v_i \rangle$, $i=1, \ldots , m+1$, and $v_{m+2} = (t_{m+2,m+1})^{-1} \widetilde{v}_{m+2}$ with $t_{m+2,m+1} = \Vert \widetilde{v}_{m+2} \Vert_2$.
Then \eqref{equ:ArnoldiTransf} can be rearranged to $K(\xi)v_{m+1} = V_{m+2} t_{m+1}$, where $V_{m+2} := [ \, V_{m+1} \; v_{m+2} \, ] \in \mathbb{R}^{dn \times (m+2)}$ and $t_{m+1} = [t_{k,m+1}]_{k=1}^{m+2} \in \RR^{m+2}$. Putting the expression $K(\xi) = ( G^2 - \xi^2 I_{dn})^{-1}$ from \eqref{equ:KInverse} in use we obtain
\begin{equation}
G^2 V_{m+2} t_{m+1} = v_{m+1} + \xi^2 V_{m+2}t_{m+1}.
\label{equ:ArnoldiTransf1} \end{equation}
The relation established in \eqref{equ:ArnoldiTransf1} can now be incorporated into the decomposition \eqref{equ:ArnoldiDec1} easily by defining
\begin{equation}
\underline{T}_{m+1}: = \left[ \begin{array}{c|c} \underline{T}_{m} & \begin{array}{c} t_{1,m+1} \\ \vdots \\ t_{m,m+1} \\ t_{m+1,m+1} \end{array} \\ \hline \begin{array}{ccc} 0 & \cdots & 0 \end{array} & t_{m+2,m+1} \end{array} \right] \in \mathbb{R}^{(m+2)\times (m+1)} \label{equ:ArnoldiExt0} \end{equation}
and
\begin{equation}
\underline{H}_{m+1} := \left[ \begin{array}{c|c} \underline{H}_{m} & \begin{array}{c} \xi^2 t_{1,m+1} \\ \vdots \\ \xi^2 t_{m,m+1} \\ 1+\xi^2 t_{m+1,m+1} \end{array} \\ \hline \begin{array}{ccc} 0 & \cdots & 0 \end{array} & \xi^2 t_{m+2,m+1} \end{array} \right] \in \mathbb{R}^{(m+2)\times (m+1)}
\label{equ:ArnoldiExt}
\end{equation}
which gives a new decomposition $G^2 V_{m+2} \underline{T}_{m+1} = V_{m+2} \underline{H}_{m+1}$
that has the same structure as in \eqref{equ:ArnoldiDec1}.

Inspired by \cite{BennEff}, we will work with a slightly modified form of the generalized rational Arnoldi decomposition from \eqref{equ:ArnoldiDec1} in all further discussions. This decomposition will turn out to be adequate for the realization of the Krylov-Schur restart discussed in Section \ref{sec:RatEvenIRA}.
To illustrate the idea, assume that $v_1 \in \RR^{dn}$ with $\Vert v_1 \Vert_2 = 1$ is given and $m=1$. The computations in \eqref{equ:ArnoldiTransf}, \eqref{equ:ArnoldiExt0}  and \eqref{equ:ArnoldiExt} yield
\begin{equation}
V_2 = \begin{bmatrix} v_1 & v_2 \end{bmatrix} \in \RR^{dn \times 2}, \quad \underline{T}_{1} = \begin{bmatrix} t_{1,1} \\ t_{2,1} \end{bmatrix} \in \RR^{2 \times 1}, \quad \underline{H}_{1} = \begin{bmatrix} h_{1,1} \\ h_{2,1} \end{bmatrix} \in \RR^{2 \times 1},
\label{equ:ArnoldiMod1}
\end{equation}
so that $G^2 V_2 \underline{T}_1 = V_2 \underline{H}_1$ holds. Now there exists a Givens rotation $F \in \RR^{2 \times 2}$ such that the second entry in $F\underline{T}_{1}$ is zero. 
Redefining $\underline{T}_1$ as $F \underline{T}_1$, $V_2$ as $V_2F^T = [ \, v_1 \; v_2 \, ]$ and $\underline{H}_{1}$ as $F\underline{H}_{1}$, we have computed a new equivalent decomposition $G^2 V_2 \underline{T}_1 = V_2 \underline{H}_1$. Now, notice that the left-hand-side can also be expressed as $G^2 V_1 T_1$, where $V_1 = [ \, v_1 \, ]$ consists only of the first column of $V_2$ and $T_1 = [ \, t_{1,1} \, ]$ where $t_{1,1}$ is the first entry of $\underline{T}_1$.
In particular, $T_1$ is now (trivially) an upper-triangular matrix.

For the further extension of the decomposition $G^2 V_1 T_1 = V_2 \underline{H}_1$, it is appropriate and feasible to keep $T_k \in \RR^{k \times k}$, $k \geq 2$, in upper-triangular form (instead of upper-Hessenberg form) throughout while the upper-Hessenberg structure of $\underline{H}_k$ is preserved. This can be achieved by applying a special bulge-chasing after every extension step. We describe this procedure in general for a given decomposition of the above form of size $m \geq 1$, i.e.
\begin{equation} \left[ \, \drawmatrix{G^2} \,  \right] \left[ \, \drawmatrix[width=.92]{V_{m}} \,  \right] \left[ \, \drawmatrix[upper]{\;\;T_{m}} \, \right] = \left[ \, \drawmatrix[width=.92]{V_{m+1}} \,  \right] \left[ \, \drawmatrix[size=1.2,width=1.1, lower banded]{\;\;\underline{H}_{m}} \, \right], \label{equ:ArnoldiMod3}
\end{equation}
where $V_{m+1} \in \RR^{dn \times (m+1)}$ has orthonormal columns $v_1, \ldots , v_{m+1}$, $T_m \in \RR^{m \times m}$ is upper-triangular and $\underline{H}_{m} \in \RR^{(m+1) \times m}$ has upper-Hessenberg structure.

Let the scalars $t_{k,m+1} \in \RR$, $k=1, \ldots , m+2$, and the vector $v_{m+2} \in \RR^{dn}$ be computed as in \eqref{equ:ArnoldiTransf} and let $\widehat{V}_{m+2} := [ \, V_{m+1} \; v_{m+2} \, ]$. The matrices in \eqref{equ:ArnoldiExt0} and \eqref{equ:ArnoldiExt} now have the special form
\begin{equation}
\widehat{\underline{T}}_{m+1} := \left[ \begin{array}{c|c} T_{m} & \begin{array}{c} t_{1,m+1} \\ \vdots \\ t_{m,m+1} \end{array} \\ \hline  \begin{array}{ccc} 0 & \cdots & 0 \\ 0 & \cdots & 0 \end{array} & \begin{array}{c} t_{m+1,m+1} \\ t_{m+2,m+1} \end{array} \end{array} \right] \in \mathbb{R}^{(m+2)\times (m+1)} \label{equ:new_TH0} \end{equation}
and
\begin{equation}
\widehat{\underline{H}}_{m+1} := \left[ \begin{array}{c|c} \underline{H}_{m} & \begin{array}{c} \xi^2 t_{1,m+1} \\ \vdots \\ \xi^2 t_{m,m+1} \\ 1+\xi^2 t_{m+1,m+1} \end{array} \\ \hline \begin{array}{ccc} 0 & \cdots & 0 \end{array} & \xi^2 t_{m+2,m+1} \end{array} \right] \in \mathbb{R}^{(m+2)\times (m+1)}.
\label{equ:new_TH}
\end{equation}
From $\widehat{V}_{m+2}$ and the upper-Hessenberg matrices in \eqref{equ:new_TH0} and \eqref{equ:new_TH} we may now recover the structures from \eqref{equ:ArnoldiMod3}, i.e. $V_{m+2} \in \RR^{dn \times (m+2)}$ with orthonormal columns, $T_{m+1} \in \RR^{(m+1) \times (m+1)}$ with upper-triangular form and $\underline{H}_{m+1} \in \RR^{(m+2) \times (m+1)}$ with upper-Hessenberg structure, so that $G^2 V_{m+1} T_{m+1} = V_{m+2} \underline{H}_{m+1}$ holds. In fact,
two orthogonal matrices $Q \in \RR^{(m+2) \times (m+2)}$ and $Z \in \RR^{(m+1) \times (m+1)}$ may be found such that $\underline{H}_{m+1} := Q\underline{\widehat{H}}_{m+1}Z \in \RR^{(m+2) \times (m+1)}$  and
\begin{equation}
Q \underline{\widehat{T}}_{m+1} Z =: \left[ \begin{array}{c} \drawmatrix[size=1.2,upper]{\;\,T_{m+1}} \\ \hline 0 \, \cdots \, 0 \end{array} \right] \quad \textnormal{with $T_{m+1} \in \RR^{(m+1) \times (m+1)}$ upper-triangular.}
\label{equ:TransformationT}
\end{equation}
Finally, defining $V_{m+2} := \widehat{V}_{m+2}Q^T$, we obtain the new decomposition
$$G^2 V_{m+1} T_{m+1} = V_{m+2} \underline{H}_{m+1}.$$
This equation is of the same form as \eqref{equ:ArnoldiMod3} except that $V_{m+1}$ has been extended by one column - which corresponds to the extension of the underlying Krylov space by one dimension - and that $T_{m+1}$ and $\underline{H}_{m+1}$ have increased in their sizes by one.
The matrices $Q$ and $Z$ can be set up as a product of Givens rotations by the bulge-chasing-process described in Algorithm \ref{alg:BC1} for a real or purely imaginary shift.

\begin{algorithm}
\caption{Bulge-Chasing-Procedure (real/imaginary shift)}
\label{alg:BC1}
\begin{algorithmic}[1]
\STATE At first, a Givens rotation is applied to $\underline{\widehat{T}}_{m+1}$ (from the left) on rows $m+2$ and $m+1$ to eliminate $t_{m+2,m}$. Applying this transformation to $\underline{\widehat{H}}_{m+1}$ introduces a bulge in the position $(m+2,m)$. This bulge can be eliminated by applying a Givens rotation (from the right) to $\underline{\widehat{H}}_{m+1}$ acting on columns $m$ and $m+1$. A bulge will now show up in the position $(m+1,m)$ in $\underline{\widehat{T}}_{m+1}$.
\STATE The bulge in the position $(m+1,m)$ in $\underline{\widehat{T}}_{m+1}$ created in (a) can be eliminated by a Givens rotation applied (from the left) on rows $m$ and $m+1$ of $\underline{\widehat{T}}_{m+1}$. This introduces a new bulge in $\underline{\widehat{H}}_{m+1}$ at the position $(m+1,m-1)$. The elimination of this bulge can be achieved by applying a Givens rotation (from the right) to $\underline{\widehat{H}}_{m+1}$ acting on the columns $m-1$ and $m$. In consequence, a new bulge will appear in $\underline{\widehat{T}}_{m+1}$ in the position $(m,m-1)$.
\STATE The elimination process described in steps 1 and 2 continues in the same manner until the bulge in $\underline{\widehat{T}}_{m+1}$ is chased off the top-left corner.
\end{algorithmic}
\end{algorithm}

Next, we discuss the case where $\xi \in \CC$ has nonzero real and imaginary parts. Let us begin directly with a real rational Arnoldi decomposition as in \eqref{equ:ArnoldiMod3}. As $K(\xi) \in \CC^{dn \times dn}$ is now a complex matrix, the resulting vector $K(\xi)v_{m+1}$ will also be complex (although $v_{m+1}$ is still real).
To remain in real arithmetics, we decompose $K(\xi)v_{m+1}$ as $\textnormal{Re}(K(\xi)v_{m+1}) + \textnormal{Im}(K(\xi)v_{m+1})i$ into its real and imaginary part. Now we apply the Gram-Schmidt process to both vectors one after the other. That is, for $\textnormal{Re}(K(\xi)v_{m+1})$ we obtain, analogously to \eqref{equ:ArnoldiTransf},
\begin{equation}
\widetilde{v}_{m+2} =\textnormal{Re}(K(\xi)v_{m+1}) - \begin{bmatrix} v_1 & \cdots & v_{m+1} \end{bmatrix} \begin{bmatrix} t_{1,m+1} \\ \vdots \\ t_{m+1,m+1} \end{bmatrix}
\label{equ:cplx_shift1}
\end{equation}
with $t_{i,m+1} = \langle  \textnormal{Re}(K(\xi)v_{m+1}),v_i \rangle$ and set $v_{m+2} := (t_{m+2,m+1})^{-1}\widetilde{v}_{m+2}$ with $t_{m+2,m+1}$ $= \Vert \widetilde{v}_{m+2} \Vert_2$. Having computed $v_{m+2}$, we may now orthogonalize $\textnormal{Im}(K(\xi)v_{m+1})$ against $v_1, \ldots , v_{m+2}$ to obtain
\begin{equation}
\widetilde{v}_{m+3} = \textnormal{Im}(K(\xi)v_{m+1}) - \begin{bmatrix} v_1 & \cdots & v_{m+2} \end{bmatrix} \begin{bmatrix} t_{1,m+2} \\ \vdots \\ t_{m+2,m+2} \end{bmatrix}
\label{equ:cplx_shift2}
\end{equation}
with $t_{i,m+2} = \langle  \textnormal{Im}(K(\xi)v_{m+1}),v_i \rangle$. Again we define $v_{m+3} := (t_{m+3,m+2})^{-1}\widetilde{v}_{m+3}$, where $t_{m+3,m+2} = \Vert \widetilde{v}_{m+3} \Vert_2$. Now we set $\widehat{V}_{m+3} = [ \, V_{m+1} \; v_{m+2} \; v_{m+3} \, ]$,

$$t_{m+1} := \begin{bmatrix} t_{1,m+1} \\ \vdots \\ t_{m+2,m+1} \\ 0 \end{bmatrix} \in \RR^{m+3} \quad \textnormal{and} \quad t_{m+2} := \begin{bmatrix} t_{1,m+2} \\ \vdots \\ t_{m+2,m+2} \\ t_{m+3,m+2} \end{bmatrix} \in \RR^{m+3}.$$
From \eqref{equ:cplx_shift1} and \eqref{equ:cplx_shift2} we obtain $ \textnormal{Re}(K(\xi)v_{m+1}) = \widehat{V}_{m+3}t_{m+1}$ and $\textnormal{Im}(K(\xi)v_{m+1}) = \widehat{V}_{m+3}t_{m+2}$, so that
$K(\xi)v_{m+1} = \widehat{V}_{m+3}(t_{m+1} + i t_{m+2})$ follows. Putting again $K(\xi) = ( G^2 - \xi^2 I_{dn})^{-1}$ from \eqref{equ:KInverse} in use we get
\begin{equation} G^2 \widehat{V}_{m+3} \big( t_{m+1} + i t_{m+2} \big) = \widehat{V}_{m+3} \big( e_{m+1} + \xi^2 t_{m+1} + i \xi^2 t_{m+2} \big), \label{equ:cplx_shift3} \end{equation}
where $e_{m+1}$ denotes the $(m+1)$-st unit vector from $\RR^{m+3}$. Furthermore, from \eqref{equ:cplx_shift3} the splitting of $\xi^2$ as $\xi^2 = \rho + \eta i$ with $\rho := \textnormal{Re}(\xi^2)$ and $\eta := \textnormal{Im}(\xi^2)$ yields
$$ G^2 \widehat{V}_{m+3} \big( t_{m+1} + i t_{m+2} \big) = \widehat{V}_{m+3} \big[ e_{m+1} + \rho t_{m+1} - \eta t_{m+2} + i \big( \eta t_{m+1} + \rho t_{m+2} \big) \big]$$
and, decomposing this once more into its real and imaginary parts, we arrive at
\begin{align}
G^2 \widehat{V}_{m+3} t_{m+1} &= \widehat{V}_{m+3} \big( e_{m+1} + \rho t_{m+1} - \eta t_{m+2} \big) \qquad \textnormal{and} \label{equ:cplx_shift4}\\
G^2 \widehat{V}_{m+3} t_{m+2} &= \widehat{V}_{m+3} \big( \eta  t_{m+1} + \rho t_{m+2} \big). \label{equ:cplx_shift5}
\end{align}
The two relations \eqref{equ:cplx_shift4} and \eqref{equ:cplx_shift5} can now be incorporated into the decomposition \eqref{equ:ArnoldiMod3}. To this end, we define
\begin{equation} \underline{\widehat{T}}_{m+2} = \left[ \begin{array}{c|c} T_{m} & \begin{array}{cc} t_{1,m+1} & t_{1,m+2} \\ \vdots & \vdots \\ t_{m,m+1} & t_{m,m+2} \end{array} \\ \hline  \begin{array}{ccc} 0 & \cdots & 0 \\ 0 & \cdots & 0 \\ 0 & \cdots & 0 \end{array} & \begin{array}{cc} t_{m+1,m+1} & t_{m+1,m+2} \\ t_{m+2,m+1} & t_{m+2,m+2} \\ 0 & t_{m+3,m+2} \end{array} \end{array} \right] \label{equ:new_THcomplex1} \end{equation}
and
\begin{equation} \underline{\widehat{H}}_{m+2} \left[ \begin{array}{c|c} \underline{H}_m & \begin{array}{cc} \rho t_{1,m+1} - \eta t_{1,m+2} & \eta t_{1,m+1} + \rho t_{1,m+2} \\ \vdots & \vdots \\ \rho t_{m,m+1} - \eta t_{m,m+2} & \eta t_{m,m+1} + \rho t_{m,m+2}  \\ 1 + \rho t_{m+1,m+1}- \eta t_{m+1,m+2} & \eta t_{m+1,m+1} + \rho t_{m+1,m+2} \end{array} \\ \hline \begin{array}{ccc} 0 & \cdots & 0 \\ 0 & \cdots & 0 \end{array} & \begin{array}{cc} \rho t_{m+2,m+1} - \eta t_{m+2,m+2} & \eta t_{m+2,m+1} + \rho t_{m+2,m+2} \\ - \eta t_{m+3,m+2} & \rho t_{m+3,m+2} \end{array} \end{array} \right]. \label{equ:new_THcomplex2} \end{equation}
From $\widehat{V}_{m+3}$ and the matrices in \eqref{equ:new_THcomplex1} and \eqref{equ:new_THcomplex2} we may again recover the structures from \eqref{equ:ArnoldiMod3}, that is
$G^2 V_{m+2} T_{m+2} = V_{m+3} \underline{H}_{m+2}$, where $V_{m+3} \in \RR^{dn \times (m+3)}$ has orthonormal columns, $T_{m+2} \in \RR^{(m+2) \times (m+2)}$ is upper-triangular and $\underline{H}_{m+2} \in \RR^{(m+3) \times (m+2)}$ has upper-Hessenberg form.
As before, a special bulge-chasing procedure is appropriate to determine two orthogonal matrices $Q \in \RR^{(m+3) \times (m+3)}$ and $Z \in \RR^{(m+2) \times (m+2)}$ such that $\underline{H}_{m+2} := Q\underline{\widehat{H}}_{m+2}Z \in \RR^{(m+3) \times (m+2)}$  and

\begin{equation}
Q \underline{\widehat{T}}_{m+2} Z = \left[ \begin{array}{c} \drawmatrix[size=1.0,upper]{\;\,T_{m+2}} \\  0 \, \cdots \, 0 \end{array} \right] \quad \textnormal{where $T_{m+2} \in \RR^{(m+2) \times (m+2)}$ is upper-triangular.}
\label{equ:TransformationT1}
\end{equation}
The matrices $Q$ and $Z$ can be set up as a product of Givens rotations by Algorithm \ref{alg:BC2}. With $V_{m+3} := \widehat{V}_{m+3}Q^T$ we obtain the desired decomposition.

\begin{algorithm}
\caption{Bulge-Chasing-Procedure (complex shift)}
\label{alg:BC2}
\begin{algorithmic}[1]
\STATE At first, a Givens rotation is applied to $\underline{\widehat{T}}_{m+2}$ (from the left) on rows $m+1$ and $m+2$ to eliminate $t_{m+2,m+1}$. Subsequently, another Givens rotation is applied to the resulting matrix on rows $m+2$ and $m+3$ to eliminate $t_{m+3,m+2}$. Applying both transformation to $\underline{\widehat{H}}_{m+2}$ introduces a bulge in the positions $(m+2,m)$ and $(m+3,m)$. 
We now apply two Givens rotations (from the right) to $\underline{\widehat{H}}_{m+2}$ acting on columns $m$ and $m+1$ to eliminate the element in position $(m+3,m)$ and, subsequently, acting on columns $m+1$ and $m+2$ to eliminate the element in position $(m+3,m+1)$. Two new bulges will show up in the positions $(m+1,m)$ and $(m+2,m+1)$ in $\underline{\widehat{T}}_{m+2}$. The additional element in the $(m+2,m)$ position of $\underline{\widehat{H}}_{m+2}$ remains in its position and is eliminated in step 2.
\STATE The bulges in the positions $(m+1,m)$ and $(m+2,m+1)$ in $\underline{\widehat{T}}_{m+2}$ created in step 1 can be eliminated by two Givens rotations applied (from the left) on rows $m$ and $m+1$ (to eliminate the bulge in position $(m+1,m)$) and on rows $m+1$ and $m+2$ (to eliminate the bulge in position $(m+2,m+1)$) of $\underline{\widehat{T}}_{m+2}$. This introduces new additional nonzero elements in $\underline{\widehat{H}}_{m+2}$ at the position $(m+1,m-1)$ and $(m+2,m-1)$. 
We apply two subsequent Givens rotations (from the right) to $\underline{\widehat{H}}_{m+2}$ acting on the columns $m-1$ and $m$ (to eliminate the element in $(m+2,m-1)$) and on columns $m$ and $m+1$ (to eliminate the element in $(m+2,m)$). Notice that the $(m+2,m)$-element we eliminate now was the one that remained in step 1. Now new bulges will appear in $\underline{\widehat{T}}_{m+2}$ in the positions $(m,m-1)$ and $(m+1,m)$. The additional element in the $(m+1,m-1)$ position of $\underline{\widehat{H}}_{m+2}$ remains in its position and is eliminated in the next step.
\STATE The elimination process described in steps 1 and 2 continues in the same manner until the bulge in $\underline{\widehat{T}}_{m+2}$ is chased off the top-left corner.
\end{algorithmic}
\end{algorithm}

Starting with some $v_1 \in \RR^{dn}, \Vert v \Vert_1 = 2,$ the previously described procedures are appropriate to construct and extend a rational Arnoldi decomposition of the form \eqref{equ:ArnoldiMod3}.
In each run, a new shift parameter $\xi \notin \sigma(\mathcal{L}_P)$ can be chosen. Per iteration, the decomposition grows in size by one if $\xi$ is real or purely imaginary and by two otherwise.

Now recall that, whenever $X$ is singular, $G$ and, consequently, a decomposition of the form \eqref{equ:ArnoldiMod3}, does not exist. Nevertheless, the vectors $K(\xi)v_{m+1}$ and $v_{m+2}$ can still be calculated as in \eqref{equ:ArnoldiTransf} and the matrices $T_{m}$, $\underline{H}_{m}$ can be extended as in \eqref{equ:new_TH0} and \eqref{equ:new_TH} if $\xi$ is real or purely imaginary. The bulge-chasing-procedure from Algorithm \ref{alg:BC1} applies and recovers the matrix structures from \eqref{equ:ArnoldiMod3}. If $\xi$ is not real or purely imaginary, $K(\xi)v_{m+1}$ can be splitted into its real and imaginary parts and the calculations in \eqref{equ:cplx_shift1} and \eqref{equ:cplx_shift2} can be carried out as described above. The extension of $T_m$ and $\underline{H}_m$ works as explained in \eqref{equ:new_THcomplex1} and \eqref{equ:new_THcomplex2} and the bulge-chasing-procedure from Algorithm \ref{alg:BC2} recovers the upper-triangular and upper-Hessenberg structures.

In conclusion, for any $m \geq 1$, the matrix pencil $\lambda T_{m} + H_{m} \in \RR^{m \times m}$ can be formed even if $G$ cannot. Moreover, its eigenvalues can be used to approximate the eigenvalues in $\sigma(\mathcal{L}_P)^2$ as before. We will permanently drop the assumption that $X$ needs to be nonsingular and that $G^2$ needs to exist from now on. In other words, we explicitly allow $\mathcal{L}_P(\lambda)$ to have eigenvalues at infinity. Therefore, the following derivations will mostly be dealing only with the matrices $V_k, \underline{T}_k$ and $\underline{H}_k$ as in \eqref{equ:ArnoldiMod3} instead of the decomposition $G^2V_mT_m = V_{m+1} \underline{H}_m$. These matrices and their modifications in the upcoming section should always be understood in the context of a rational Arnoldi decomposition as in \eqref{equ:ArnoldiDec0} whenever such a decomposition exists.
We have summarized the method to generate (or extend) a rational Arnoldi decomposition in Algorithm \ref{alg:RatArnoldi}.
\begin{algorithm}
\caption{Rational Arnoldi Expansion}
\label{alg:RatArnoldi}
\begin{algorithmic}[1]
\STATE \textsc{Input}: The linearization $\mathcal{L}_P(\lambda) = \lambda X + Y \in \RR[\lambda]^{dn \times dn}$ for a $T$-even matrix polynomial $P(\lambda) \in \RR[\lambda]^{n \times n}$ defined in \eqref{equ:Linearization}. A matrix $V_{k+1} =[ \, v_1 \; \cdots \; v_{k+1} \, ] \in \RR^{dn \times (k+1)}$ with orthonormal columns, $T_k \in \RR^{k \times k}$ upper-triangular and $\underline{H}_k \in \RR^{(k+1) \times k}$ in upper Hessenberg form satisfying $G^2V_kT_k = V_{k+1} \underline{H}_k$ if $G^2$ exists. In case $k=0$, we set $T_0 = [ \, ]$ and $\underline{H}_0 := [ \, ]$. A number $m \in \mathbb{N}$, $m > k$.
\STATE \textsc{Output}: Matrices $V_{m+1} =[ \, v_1 \; \cdots \; v_{m+1} \, ] \in \RR^{dn \times (m+1)}$ with orthonormal columns, $T_m \in \RR^{m \times m}$ upper-triangular and $\underline{H}_m \in \RR^{(m+1) \times m}$ in upper Hessenberg form that satisfy \eqref{equ:ArnoldiDec1} in case $G^2$ exists.
       \FOR {$j = \ell+1, \ldots , m$}
       \STATE pick a shift $\zeta_j  \in \CC$
       \STATE compute $w := K(\zeta_j)v_j = (\mathcal{L}_P(\zeta_j)^{-T}X \mathcal{L}_P(\zeta_j)^{-1}X)v_j$ using Section \ref{sec:LinearSystems}
       \IF {$\zeta_j \in \RR$ or $\zeta_j \in i \RR$}
       \STATE orthogonalize $w$ against $v_1, \ldots , v_j$  and obtain $t_{1,j}, \ldots , t_{j,j} \in \RR$ as in \eqref{equ:ArnoldiTransf}
       \STATE set $v_{j+1}$ to obtain $t_{j+1,j} \in \RR$
       \STATE form $\underline{\widehat{T}}_j \in \RR^{(j+1) \times j}$ and $\underline{\widehat{H}}_j \in \RR^{(j+1) \times j}$ as in \eqref{equ:new_TH0} and \eqref{equ:new_TH}
       \STATE set $\widehat{V}_{j+1} = [ \, v_1 \; \cdots \; v_{j+1} \, ]$
       \STATE apply the bulge-chasing-procedure described in Algorithm \ref{alg:BC1} to determine
       \STATE orthogonal matrices $Q \in \RR^{(j+1) \times (j+1)}$ and $Z \in \RR^{j \times j}$ such that
       \STATE \hspace{1cm} $\triangleright$ $Q\underline{\widehat{T}}_jZ$ is upper-triangular with zeros in its last row and
       \STATE \hspace{1cm} $\triangleright$ $Q\underline{\widehat{H}}_jZ =: \underline{H}_j$ has upper-Hessenberg structure
       \STATE define $T_j$ to be the first $j$ rows of $Q\underline{\widehat{T}}_jZ$
       \STATE define $V_{j+1} := \widehat{V}_{j+1}Q^T$
       \ELSE
       \STATE orthogonalize $\textnormal{Re}(w)$ against $v_1, \ldots , v_j$  and get $t_{1,j}, \ldots , t_{j,j} \in \RR$ as in \eqref{equ:cplx_shift1}
       \STATE set $v_{j+1}$ to obtain $t_{j+1,j} \in \RR$
       \STATE orthogonalize $\textnormal{Im}(w)$ against $v_1, \ldots , v_j, v_{j+1}$  and obtain $t_{1,j+1}, \ldots , t_{j+1,j+1} \in \RR$ as  in \eqref{equ:cplx_shift2}
       \STATE set $v_{j+2}$ to obtain $t_{j+2,j+1} \in \RR$
       \STATE form $\underline{\widehat{T}}_{j+1} \in \RR^{(j+2) \times (j+1)}$ and $\underline{\widehat{H}}_j \in \RR^{(j+2) \times (j+1)}$ as in \eqref{equ:new_THcomplex1} and \eqref{equ:new_THcomplex2}
       \STATE set $\widehat{V}_{j+2} = [ \, v_1 \; \cdots \; v_{j+1} \; v_{j+2} \, ]$
       \STATE apply the bulge-chasing-procedure described in Algorithm \ref{alg:BC2} to determine
       \STATE orthogonal matrices $Q \in \RR^{(j+2) \times (j+2)}$ and $Z \in \RR^{(j+1) \times (j+1)}$ such that
       \STATE \hspace{1cm} $\triangleright$ $Q\underline{\widehat{T}}_{j+1}Z$ is upper-triangular with zeros in its last row
       \STATE \hspace{1cm} $\triangleright$ $Q\underline{\widehat{H}}_{j+1}Z =: \underline{H}_{j+1}$ has upper-Hessenberg structure
       \STATE define $T_{j+1}$ to be the first $j+1$ rows of $Q\underline{\widehat{T}}_{j+1}Z$
       \STATE define $V_{j+2} := \widehat{V}_{j+2}Q^T$
       \ENDIF
       \ENDFOR
\end{algorithmic}
\end{algorithm}
\section{The rational Even-IRA algorithm}
\label{sec:RatEvenIRA}
Let $P(\lambda) \in \RR[\lambda]^{n \times n}$ be some $T$-even matrix polynomial and let $\mathcal{L}_P(\lambda) = \lambda X + Y \in \RR[\lambda]^{dn \times dn}$ and
$$ K(\zeta) = \mathcal{L}_P(\zeta)^{-T}X \mathcal{L}_P(\zeta)^{-1}X = \mathcal{L}_P(-\zeta)^{-1}X \mathcal{L}_P(\zeta)^{-1}X, \quad \zeta \notin \sigma(P),$$
be defined for $P(\lambda)$ as in \eqref{equ:Linearization} and \eqref{equ:SpecTransform}, respectively. The rational \textsc{Even-IRA} algorithm presented in this section is a method that unifies the Krylov-Schur restart strategy  \cite{KrylovSchur} with the spectral-preserving transformation $K(\zeta)$ (see Section \ref{sec:EvenIRA} and \cite{EvenIRA,MehrWat01,Wat04}) and the shift flexibility offered by the rational Arnoldi process \cite{Ruhe98, BerlGuett15, Ruhe84}.  The sparse and structured form of the linearization $\mathcal{L}_P(\lambda)$ (see Theorem \ref{thm:linearization}) is exploited for evaluating matrix-vector-products with $K(\zeta)$ implicitly and efficiently without ever forming $K(\zeta)$ at all (see Section \ref{sec:LinearSystems}). Hence, the memory requirement of the method is essentially that of storing the given matrix polynomial and the vectors from the current search space. In a nutshell, this approach yields a powerful Krylov-subspace algorithm for the computation of some eigenvalues for $T$-even polynomial eigenvalue problems. The rational Even-IRA algorithm presented next consists of several phases. In the initialization phase (Section \ref{ssec:initialization}) a rational Arnoldi decomposition is constructed which is the start and end point of each Krylov-Schur cycle. In the expansion phase (Section \ref{ssec:expansion}) the size of this decomposition is increased. After the expansion, a QZ decomposition is applied (Section \ref{ssec:decomp}) to identify eigenvalues that have converged during the current run and which are to be locked (Section \ref{ssec:convergence}). To initialize the algorithm's next cycle, the decomposition is truncated (Section \ref{ssec:truncation}) and the upper-triangular and upper-Hessenberg forms of the matrices are recovered (Section \ref{ssec:recovery}). The next iteration then begins with the expansion phase. We now describe the different phases in detail.

\subsection{The initialization phase}
\label{ssec:initialization}
Let $M \in \mathbb{N}$ be the number of desired eigenvalues for $P(\lambda)$ ($\mathcal{L}_P(\lambda)$, respectively).  The premier start of the algorithm begins with its initialization phase. That is, matrices
\begin{equation}
 \left[ \, \drawmatrix[width=1,height=1.2]{V_{M+1}} \,  \right] \in \RR^{dn \times (M+1)}, \; \left[ \, \drawmatrix[upper, size=1.2]{\quad T_M} \, \right] \in \RR^{M \times M} , \;   \underline{H}_M = \left[ \begin{array}{c}  H_M  \\[0.1cm] \hline B_M \end{array} \right] \in \RR^{(M+1) \times M},
\label{equ:Restart-1}
\end{equation}
where $B= h_{M+1,M} e_M^T$ for some scalar $h_{M+1,M} \in \RR$ are computed by Algorithm \ref{alg:RatArnoldi}. Note that the columns of $V_{m+1}$ are orthonormal, $T_M$ is upper-triangular and $H_M$ has upper-Hessenberg form. In case $X$ is nonsingular, $(X^{-1}Y)^2 V_MT_M = V_{M+1}\underline{H}_M$ holds. If $P(\lambda)$ has no eigenvalues at infinity, Algorithm \ref{alg:RatArnoldi} may be initialized with $V_1 = [\, v_1 \,]$ (arbitrary and normalized), $T_0 = [ \, ]$ and $\underline{H}_0 = [ \, ]$. In case of the presence of infinite eigenvalues, a different initialization should be chosen, see Section \ref{sec:EigenvalueInfinity}.

A cycle of the rational Krylov-Schur algorithm begins and ends with matrices of the form \eqref{equ:Restart-1}. Now suppose, at some stage of the algorithm, $s \in \mathbb{N}_0$ eigenvalues have already converged. Assume these had been locked so that they are located in the top-left $s \times s$ corner of $T_M$ and $H_M$ (of course, beginning with the algorithms first run, $s=0$).

\subsection{The expansion phase}
\label{ssec:expansion}
The first step of the algorithm is the expansion phase where the above matrices are extended up to a size $m > M$. This is achieved by performing $m-M$ additional steps of Algorithm \ref{alg:RatArnoldi} with the input matrices from \eqref{equ:Restart-1}. We call $m - M$ the extension size for the algorithm. Now we obtain matrices
\begin{equation}
 \left[ \, \drawmatrix[width=1,height=1.2]{V_{m+1}} \,  \right] \in \RR^{dn \times (m+1)}, \; \left[ \, \drawmatrix[upper, size=1.2]{\quad T_m} \, \right] \in \RR^{m \times m} , \; \underline{H}_m = \left[ \begin{array}{c}  H_m  \\[0.1cm] \hline B_m \end{array} \right] \in \RR^{(m+1) \times m},
\label{equ:Restart0} \end{equation}
where $B_m = h_{m+1,m}e_m^T$, the columns $v_1, \ldots , v_{m+1} \in \RR^{dn}$ of $V_{m+1} = [ \, V_m \; v_{m+1} \, ]$ are orthonormal, $T_m$ is upper-triangular and $H_m$ has upper-Hessenberg structure.
We partition the matrices in \eqref{equ:Restart0} in accordance with the number $s$ of locked eigenvalues as
\begin{equation}
V_m = \begin{bmatrix} V_s & V \end{bmatrix}, \quad T_m = \begin{bmatrix} T_s & T' \\ 0 & T \end{bmatrix}, \quad H_m = \begin{bmatrix} H_s & H' \\ 0 & H \end{bmatrix}  \quad \textnormal{and} \quad B_m = \begin{bmatrix} B_s & B \end{bmatrix},
\label{equ:Restart1}
\end{equation}
where $T_s, H_s \in \RR^{s \times s}$, $B_s^T = [ \, 0 \; \cdots \; 0 \, ]^T \in \RR^s$ and $B^T = [ \, 0 \; \cdots \; 0 \; h_{m+1,m} \, ]^T \in \RR^k$, where we have set $k := m-s$ implying $T,H \in \RR^{k \times k}$. Moreover, recall that the eigenvalues from the matrix pair $\lambda T_s + H_s$ are the ones we assumed to be locked.

\subsection{The decomposition and reordering phase}
\label{ssec:decomp}
We may now enter the decomposition and reordering phase of the algorithm. To this end, we first compute a QZ decomposition of the matrix pair $(T, H)$. For this purpose, orthogonal matrices $Q_1, Z_1 \in \RR^{k \times k}$ can be determined so that $Q_1^TTZ_1 = T^\star \in \RR^{k \times k}$ remains upper-triangular while $Q_1^THZ_1 = H^\star \in \RR^{k \times k}$ becomes quasi upper-triangular (with solely $1 \times 1$ and $2 \times 2$ blocks along its diagonal). At this point, a reordering procedure (see, e.g., \cite{Reorder}) can be applied to $\lambda T^\star + H^\star$ to move unwanted eigenvalues of $\lambda T^\star + H^\star$ into the trailing part of its generalized Schur decomposition.
That is, two additional orthogonal transformations $Q_2, Z_2 \in \RR^{k \times k}$ can be found, so that unwanted eigenvalues of $\lambda T^\star + H^\star$ move to the south-east corner of the matrices $T^{\diamond} := Q_2^TT^\star Z_2 \in \RR^{k \times k}$ and $H^{\diamond} := Q_2^TH^\star Z_2 \in \RR^{k \times k}$. Thereby, the matrices $T^\diamond \in \RR^{k \times k}$ and $H^\diamond \in \RR^{k \times k}$ stay upper-triangular and quasi upper-triangular, respectively. Finally, defining $Q^T := Q_2^TQ_1^T$ and $Z := Z_1Z_2$,  we update \eqref{equ:Restart1} as follows
\begin{equation}
\widehat{V}_m = \begin{bmatrix} V_s & VQ \end{bmatrix}, \quad \widehat{T}_m = \begin{bmatrix} T_s & T'Z \\ 0 & T^\diamond \end{bmatrix}, \quad  \widehat{H}_m = \begin{bmatrix} H_s & H'Z \\ 0 & H^\diamond \end{bmatrix} \quad \textnormal{and} \quad
\widehat{B}_m = \begin{bmatrix} B_s & B^\diamond \end{bmatrix}
\label{equ:Restart2}
\end{equation}
with $B^\diamond := B^TZ$. Notice that $B^\diamond$ will now be, in general, a full vector.

\subsection{The inspection-of-convergence phase}
\label{ssec:convergence}
With \eqref{equ:Restart2} the inspection-of-convergence phase of the algorithm begins. That is, the leading components of $B^\diamond$ are inspected for convergence and eigenvalues  are locked whenever convergence has taken place. Let $H^\diamond = [ h_{i,j} ]_{i,j}, T^\diamond = [ t_{i,j} ]_{i,j}$ with $1 \leq i,j \leq k $ and let $B^\diamond = [ \, b_{s+1} \; \cdots \; b_m \, ]$. Starting with $r \equiv 1$ we now consider the following cases: \vspace*{0.1cm}
\begin{compactenum}[(a)]
\item Whenever $h_{r+1,r} = 0$ and $| b_{s+r} |$ is below a given tolerance \texttt{tol}, we consider the corresponding eigenvalue $h_{r,r}/t_{r,r}$ as converged. The element $b_{s+r}$ is set to zero and the number $r$ of converged eigenvalues in the current run is increased by one.
\item Whenever $h_{r+1,r} \neq 0$ but $\Vert [ \, b_{s+r} \; b_{s+r+1} \, ] \Vert_2$ is below the given tolerance \texttt{tol}, we consider the pair of complex conjugate eigenvalues corresponding to the $2 \times 2$ matrix pencil
    $$ \lambda \begin{bmatrix} t_{r,r} & t_{r,r+1} \\ 0 & t_{r+1,r+1} \end{bmatrix} + \begin{bmatrix} h_{r,r} & h_{r,r+1} \\ h_{r+1,r} & h_{r+1,r+1} \end{bmatrix}$$
    as converged. The elements $b_{s+r}$ and $b_{s+r+1}$ are both set to zero. Finally, the number $r$ of converged eigenvalues in the current run is increased by two.
\end{compactenum}
\vspace*{0.1cm}We repeat the locking of eigenvalues as long as (a) or (b) reveals convergence. Once no further convergence is observed notice that $\widehat{B}_m := [ \, 0 \; \cdots \; 0 \; b_{s+r+1} \; \cdots \; b_m \, ]$ (where $r$ is now the total number of locked eigenvalues during the current run). The new number of converged eigenvalues in total is now $s^\star = s+ r$. If $s^\star \geq M$ (the number of desired eigenvalues), we are done. Otherwise the matrices are truncated to prepare a restart of the algorithm.

\subsection{The truncation phase}
\label{ssec:truncation}
If $s^\star < M$, the size of the matrices in \eqref{equ:Restart2} is now decreased to size $M \times M$ in the truncation phase to initialize a restart of the process.
In particular, let $\widehat{V}_{M} \in \RR^{dn \times M}$ be the first $M$ columns of $\widehat{V}_m$ and $\widehat{V}_{M+1} = [ \, \widehat{V}_M \; v_{m+1} \, ]$, where $v_{m+1}$ denotes the last column from $V_{m+1}$ in \eqref{equ:Restart0} (note that $v_{m+1}$ has not been touched in all steps up to this point). Moreover, denote the top-left $M \times M$ submatrices of $\widehat{T}_m$ and $\widehat{H}_m$ by $\widehat{T}_M$ and $\widehat{H}_M$, respectively, and the vector obtained from the first $M$ components of $\widehat{B}_m$ by $\widehat{B}_M$, i.e. $\widehat{B}_M = [ \, 0 \; \cdots \; 0 \; b_{s^\star+1} \; \cdots \; b_M \, ]$. In the form \eqref{equ:Restart0} we have
\begin{equation}
 \left[ \, \drawmatrix[width=1,height=1.2]{\widehat{V}_{M+1}} \,  \right] \in \RR^{dn \times (M+1)}, \; \left[ \, \drawmatrix[upper, size=1.2]{\quad \widehat{T}_M} \, \right] \in \RR^{M \times M} , \; \widehat{\underline{H}}_M := \left[ \begin{array}{c}  \widehat{H}_M  \\[0.1cm] \hline \widehat{B}_M \end{array} \right] \in \RR^{(M+1) \times M}.
\notag \end{equation}
Notice that $\widehat{\underline{H}}_M$ will not have Hessenberg structure at this stage of the algorithm since $\widehat{B}_M$ will have more nonzero elements than just $b_M$.

\begin{remark}
It is unfortunate to truncate the matrices as above whenever the $(M+1,M)$-element in $\widehat{H}_m$ is nonzero. In this case a $2 \times 2$ block is split which should be avoided by decreasing or increasing $M$ by one. \end{remark}

Analogously to \eqref{equ:Restart1} and \eqref{equ:Restart2} we may now partition $\widehat{V}_{M}, \widehat{T}_M$ and $\widehat{\underline{H}}_M$ according to the new number $s^\star$ of locked and converged Ritz values. This highlights the active part of the decomposition and separates it from the locked part (which does not need to be touched again). In particular, we partition
\begin{equation}
\widehat{V}_{M} = \begin{bmatrix} V_{s^\star} & V^\circ \end{bmatrix}, \quad \widehat{T}_{M} = \begin{bmatrix} T_{s^\star} & T'' \\ 0 & T^\circ \end{bmatrix}, \quad  \widehat{H}_M = \begin{bmatrix} H_{s^\star} & H'' \\ 0 & H^\circ \end{bmatrix}, \quad  \widehat{B}_M = \begin{bmatrix} B_{s^\star} & B^\circ \end{bmatrix},
\label{equ:Restart4}
\end{equation}
where $V_{s^\star} \in \RR^{dn \times s^\star}, T_{s^\star}, H_{s^\star} \in \RR^{s^\star \times s^\star}$ and $B_{s^\star}^T = [0 \; \cdots \; 0 \, ]^T \in \RR^{s^\star}$. Recall that $B^\circ$ is in general a full vector with all nonzero entries. Set $k^\star = M - s^\star$ so that $T^\circ, H^\circ \in \RR^{k^\star \times k^\star}$.

\subsection{The recovery phase}
\label{ssec:recovery}
Our next goal is to tranform the matrices in \eqref{equ:Restart4} back to a decomposition of the form \eqref{equ:Restart-1} in the recovery phase.
That is, we determine orthogonal matrices $Q,Z \in \RR^{k^\star \times k^\star}$ such that $Q^TT^\circ Z =: T \in \RR^{k^\star \times k^\star}$ is still upper-triangular, $Q^TH^\circ Z =:H \in \RR^{k^\star \times k^\star}$ remains in upper-Hessenberg form and $B = B^\circ Z =: h_{M+1,M}e_{k^\star}^T$ for some scalar $h_{M+1,M} \in \RR$. Then we update \eqref{equ:Restart4} to obtain
\begin{equation}
V_{M} = \begin{bmatrix} V_{s^\star} & V^\circ Q \end{bmatrix}, \quad T_{M} = \begin{bmatrix} T_{s^\star} & T''Z \\ 0 & T \end{bmatrix}, \quad  H_M = \begin{bmatrix} H_{s^\star} & H''Z \\ 0 & H \end{bmatrix}, \;\,
B_M = \begin{bmatrix} B_{s^\star} & B \end{bmatrix}
\end{equation}
and with $V_{M+1} = [ \, V_M \; v_{m+1} \, ]$ we are back with matrices
\begin{equation}
 \left[ \, \drawmatrix[width=1,height=1.2]{V_{M+1}} \,  \right] \in \RR^{dn \times (M+1)}, \; \left[ \, \drawmatrix[upper, size=1.2]{\quad T_M} \, \right] \in \RR^{M \times M} , \;  \underline{H}_M = \left[ \begin{array}{c}  H_M  \\[0.1cm] \hline B_M \end{array} \right] \in \RR^{(M+1) \times M}
\label{equ:Restart5}
 \end{equation}
as in \eqref{equ:Restart-1}, where $B_M = h_{M+1,M}e_M^T$. The next cycle of the algorithm then begins with the expansion phase as described in Section \ref{ssec:expansion}. The recovery phase can be carried out by the bulge-chasing-process described in Section \ref{sec:RecoveringPhase}

The overall goal of this algorithm is to achieve $B_M=[ \,  0 \; \cdots \; 0 \, ]$ in \eqref{equ:Restart5} after some cycles of the restarting procedure described above. As soon as this situation takes place, the $M$ eigenvalues of $\lambda T_M + H_M$ are exact eigenvalues of $\mathcal{L}_P(\lambda)^2$ and, in turn, their plus/minus square roots exact eigenvalues of $P(\lambda)$.

\subsection{The eigenvalue infinity}
\label{sec:EigenvalueInfinity}

A matrix polynomial $P(\lambda)$ might have eigenvalues at infinity (see Section \ref{sec:Notation}).  The rational \textsc{Even-IRA} algorithm will eventually detect infinite eigenvalues, i.e., in computations in real arithmetic, eigenvalues of very large magnitude might be found. This is detrimental for the algorithm's performance since $(i)$ the detection of very large eigenvalues is, in this case, a wrong result, and $(ii)$ the convergence results after the detection of such an eigenvalue are of unsatisfying accuracy. Therefore, it seems reasonable to a priori eliminate any possibility of convergence to infinity. This will guarantee a good performance throughout and reliable results.

Assuming $P(\lambda) = \sum_{k=1}^d P_k \lambda^k \in \RR[\lambda]^{n \times n}$ of degree $d \geq 1$ is regular, the eigenvectors for the eigenvalue $\mu = \infty$ are the nullvectors of $P_d$. These can be found by solving the $n \times n$ linear system $P_dx=0$ with an appropriate method. These vectors can now be used to initialize our algorithm so that convergence for the eigenvalue infinity has already taken place. For this purpose, let $\dim ( \textnormal{null}(P_d)) = t$ and $\lbrace v_1, v_2, \ldots , v_t \rbrace \subset \RR^n$ some orthonormal basis of $\textnormal{null}(P_d)$.  
We define
\begin{equation}
V_{t+1} = \left[ \begin{array}{c|c} \begin{array}{ccc} v_1 & \cdots & v_t \\ 0 & & 0 \\ \vdots & & \vdots \\ 0 & \cdots & 0 \end{array} & v_{t+1} \end{array} \right] \in \RR^{dn \times (t+1)}, \quad \underline{H}_t = \left[ \begin{array}{c}  I_t  \\[0.1cm] \hline \begin{array}{ccc} 0 & \cdots & 0 \end{array} \end{array} \right],
\notag
\end{equation}
and $T_{t} = 0_{t \times t} \in \RR^{t \times t}$. The vector $v_{t+1} \in \RR^{dn}$ can be chosen arbitrarily so that the columns of $V_{t+1}$ are orthonormal. From here on, we start the Initialization Phase of the rational \textsc{Even-IRA} algorithm described in Section \ref{sec:RatEvenIRA} with $V_{t+1}$, $T_t$ and $\underline{H}_t$ in Algorithm \ref{alg:RatArnoldi}. Moreover, we define, right from this point on, the number $s$ of  converged eigenvalues to be $t$. In other words, with this initialization, convergence to infinity and locking has already occurred before the algorithm actually starts. The algorithm will not reveal further eigenvalues at infinity. The number of desired eigenvalues has to be increased from $M$ to $t+M$.

\subsection{The shift-strategy} \label{sec:shift}
The appropriate choice of shifts is a delicate problem that often depends on user-specified priorities. According to Section \ref{sec:LinearSystems}, a matrix-vector-multiplication with $K(\zeta) \in \CC^{dn \times dn}$ essentially reduces to a system solve with $P(\zeta)$ (and $P(\zeta)^T$). If an LU decomposition of $P(\zeta)$ is computed, it can be reused as long as the shift does not change. On the other hand, once a change of shift took place, a new decomposition has to be computed for the subsequent iterations (until the next change). Hence, in order to keep the algorithm effective, shift changes should not be applied too often. However, on the other hand, the choice of a good new shift is likely to increase the convergence speed. Two general shift strategies are given below.
\begin{enumerate}[(a)]
\item Assume the current run of the rational \textsc{Even-IRA} algorithm revealed convergence and (in total) $s^\star$ eigenvalues are locked - this corresponds to the situation \eqref{equ:Restart4}. Let $T^\circ = [t^\circ_{i,j}]_{i,j}, H^\circ = [h^\circ_{i,j}]_{i,j}$, $B^\circ = [ \, b_{s^\star +1} \; \cdots \; b_M \, ]$ and consider the case $h^\circ_{2,1} = 0$. In particular, assuming $G^2$ exists and regarding $\widehat{V}_{s^\star +1}$ (the matrix consisting of the first $s^\star +1$ columns of $\widehat{V}_M$), $\widehat{T}_{s^\star +1}, \widehat{H}_{s^\star +1}$ (the $(s^\star +1) \times (s^\star +1)$ principal submatrices of $\widehat{T}_M$ and $\widehat{H}_M$, respectively) we have in view of \eqref{equ:Restart4}
    $$G^2 V_{s^\star +1} T_{s^\star +1} = V_{s^\star +1} H_{s^\star +1} + b_{s^\star +1} v_{m+1} e_{s^\star +1}^T.$$
    In other words,
    \begin{equation} \Vert G^2 V_{s^\star +1} T_{s^\star +1} - V_{s^\star +1} H_{s^\star +1} \Vert_2 =  | b_{s^\star +1} | \label{equ:residual} \end{equation}
    because $\Vert v_{m+1} \Vert_2 = 1$. Therefore, the absolute value of the first entry $b_{s^\star +1}$ of $B^\circ$ displays the first residual which was not below the given tolerance \texttt{tol} since, otherwise, the corresponding eigenvalue $\xi := h^\circ_{1,1}/t^\circ_{1,1}$ of $\lambda T^\circ + H^\circ$ located in the top-left $1 \times 1$ block would have been identified as converged.
    Nevertheless, $\xi$ may serve as a good approximation to the next eigenvalue that is about to converge and $\xi$ might now be chosen as the next shift parameter. Analogously, whenever $h^\circ_{2,1} \neq 0$, an eigenvalue of the $2 \times 2$ top-left corner of $\lambda T^\circ + H^\circ$ can be chosen as a new shift.
\item The shift-strategy from (a) can be modified so that the next shift parameter is chosen as $\xi = h^\circ_{1,1}/t^\circ_{1,1}$ only if the corresponding residual $| b_{s^\star +1} |$ is above a given tolerance. In particular, if $| b_{s^\star +1} |$ is already very small, a change of shift is probably not necessary since the algorithm seems to be ``on the right way'' to reveal the next convergence soon (e.g. within the next cycle). However, $| b_{s^\star +1} |$ being above a given tolerance might indicate that the current shift is not heading off to reveal further convergence in the near future. Thus, changing the shift could be an appropriate means to speed up convergence in such a situation.
\end{enumerate}
Certainly, other shift strategies beside (a) and (b) above and mixtures of both are conceivable. In particular, if one is interested in eigenvalues in a particular region of the complex plane, the shift should, of course, be chosen appropriately.

\subsection{The recovery phase}
\label{sec:RecoveringPhase}

We now consider the recovery phase of the rational \textsc{Even-IRA} algorithm in detail. Therefore, reconsider the matrices obtained in \eqref{equ:Restart4}.
We now show how to construct two orthogonal matrices $Q,Z \in \RR^{k^\star \times k^\star}$ such that $Q^TT^\circ Q = T \in \RR^{k^\star \times k^\star}$ remains upper-triangular, $Q^TH^\circ Z = H \in \RR^{k^\star \times k^\star}$ has upper-Hessenberg form and $B^\circ Z = [ \, 0 \; \cdots \; 0 \; h_{M+1,M} \, ]$ is a vector of zeros except for some scalar $h_{M+1,M} \in \RR$ in the last position. The matrices $Q$ and $Z$ are the products of a sequence of Givens rotations that constitute our bulge-chasing procedure. Hereby, a Givens rotation $\tilde{Q}^T \in \RR^{2 \times 2}$ from the left acts on two rows $i$ and $j$ (with $1 \leq i,j \leq k^\star$) of $T^\circ$ and $H^\circ$. Each transformation $\tilde{Q}^T$ needs to be applied via $\tilde{Q}$ to the columns $i$ and $j$ of $V^\circ$, too. This is implicitly understood in all the following derivations. A Givens rotation $\tilde{Z} \in \RR^{2 \times 2}$ from the right acts on two columns $i$ and $j$, $1 \leq i,j \leq k^\star$, of $T^\circ$, $H^\circ$ and $B^\circ$. These transformations do not influence the matrix $V^\circ$.

Now let $B^\circ = [ \, b_1 \; \cdots \; b_{k^\star} \, ]$. The bulge-chasing process proceeds as follows: \vspace*{0.1cm}
\begin{compactenum}[(a)]
\item We apply a Givens rotation $Z_1$ from the right on the first two columns to eliminate $b_1$ using $b_{2}$. This introduces a bulge in the position $(3,1)$ in $H^\circ$ and in the position $(2,1)$ in $T^\circ$. A rotation $Q_1^T$ acting on rows one and two from the left can be used to eliminate the bulge in $T^\circ$. The bulge on the second subdiagonal in $H^\circ$ remains in its position. Now the first element in $B^\circ$ is zero and the analogous process can be used to eliminate the second element in $B^\circ$. As before, the new bulge in the position $(4,2)$ in $H^\circ$ (i.e. on the second subdiagonal in $H^\circ$) remains it its position.
\item The third element in $B^\circ$ can be eliminated as in (a) above and two new elements in the positions $(4,3)$ in $T^\circ$ and $(5,3)$ in $H^\circ$ show up. As in (a), this bulge in $H^\circ$ is accepted for the moment. However, with the elimination of the bulge in $T^\circ$ with a Givens rotations from the left on rows three and four an additional bulge in $H^\circ$ will appear in the position $(4,1)$ (i.e. on the third subdiagonal in $H^\circ$). This bulge can be eliminated by a Givens rotation from the right on the first and second column of $H^\circ$  introducing again a bulge in the $(2,1)$ position in $T^\circ$. A rotation applied to the first two rows from the left is used to eliminate the bulge in $T^\circ$.
\item The process from (b) now continues for all $t > 3$. That is, the elimination of the $t$-th entry in $B^\circ$ is achieved by a Givens rotation from the right on columns $t$ and $t+1$. Consequently, bulges appear in $(t+2,t)$ in $H^\circ$ and $(t+1,t)$ in $T^\circ$. The elimination of the bulge in $T^\circ$ by a Givens rotation from the left introduces an additional bulge $(t+1,t-2)$ in $H^\circ$. This bulge is chased off the top-left corner of $T^\circ$ and $H^\circ$ by applying Givens rotation alternatingly from left and right. \vspace*{0.1cm}
\end{compactenum}
If the bulge-chasing process described in (a) to (c) is completely carried out, in the end, $T^\circ$ is still of upper-triangular form, $B^\circ = [ \, 0 \; \cdots \; 0 \; h_{M+1,M} \, ]$ and $H^\circ$ is a matrix that now has two full subdiagonals (i.e. all entries below the second subdiagonal of $H^\circ$ are zero). Now the transformation process can be continued and the second subdiagonal in $H^\circ$ can be eliminated from the lower right corner to the top-left corner. A standard bulge-chasing (Givens rotations alternatingly from left and right) is adequate to achieve this. It is important to note that no Givens rotation is required that touches the last column of $T^\circ$, $H^\circ$ and $B^\circ$. Therefore, $B^\circ$ remains as it is and we obtain the desired form.

\section{Numerical experiments}
\label{sec:examples}
In this section, we briefly describe the results of two numerical experiments to give a proof of concept for the algorithm described in the previous section. To this end, we set up a basic implementation of the rational \textsc{Even-IRA} algorithm in MATLAB R2020a and compared our results to those found with the MATLAB function \texttt{polyeig}. As the degree of $P(\lambda)$ is even in both examples, $\mathcal{L}_P(\lambda)$ was constructed as in \eqref{equ:Linearization} with $M_P(\lambda)$ from Definition \ref{def:MP} (b). We initialize the algorithm as explained in Section \ref{sec:EigenvalueInfinity}. In contrast to the computation of eigenvalues with \texttt{polyeig}, the rational \textsc{Even-IRA} algorithm is designed to find only a few eigenvalues of a matrix polynomial. Therefore, a comparison of the computational times for both algorithms seems inappropriate here.

\begin{figure}
  \centering
  \includegraphics[width=170pt]{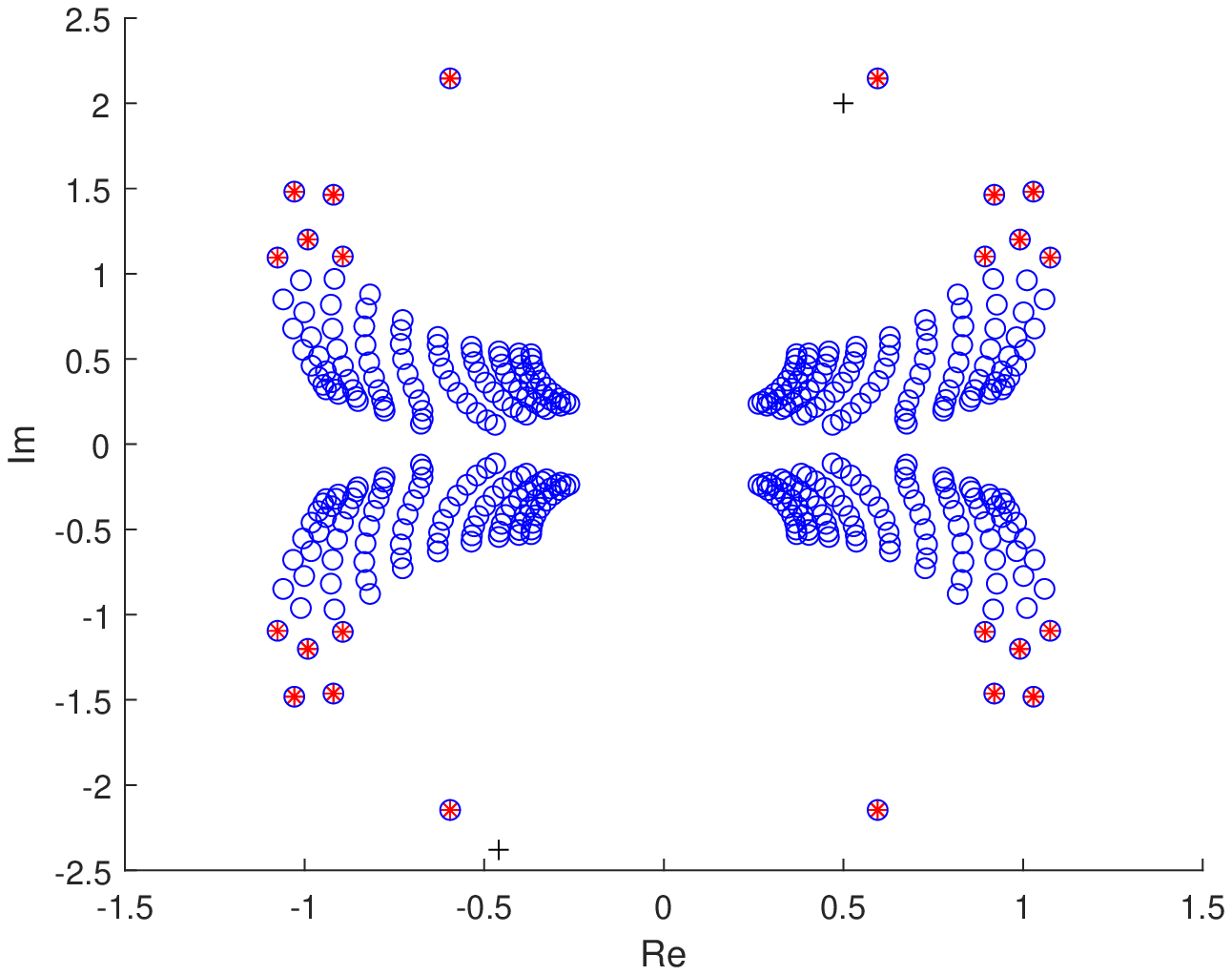} \; \includegraphics[width=170pt]{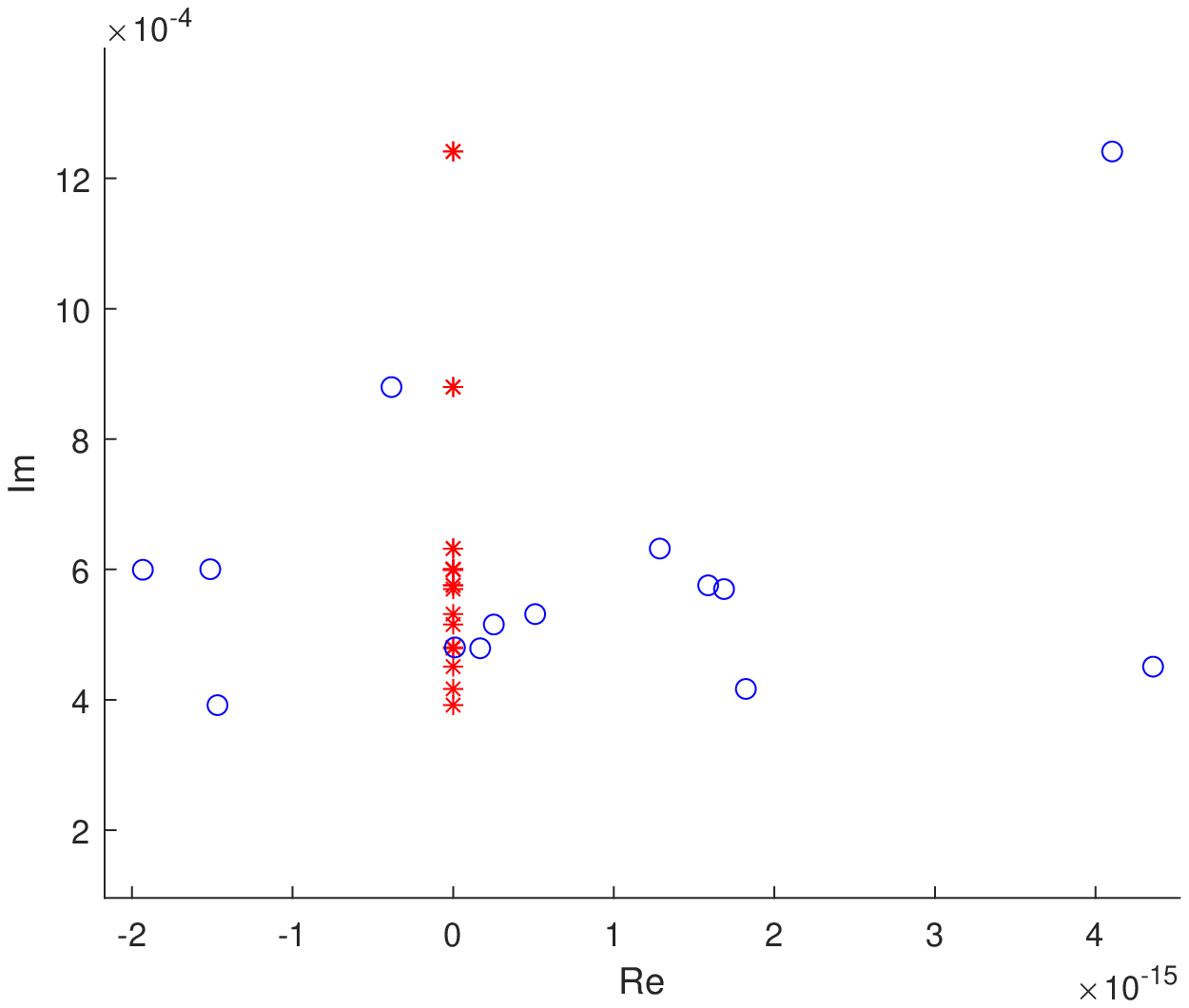}
  \caption{Left: 24 eigenvalues found by the rational \textsc{Even-IRA} algorithm (red stars) for the example \texttt{butterfly} in \cite{Bet13}. Blue circles indicate the eigenvalues computed via \texttt{polyeig}. Due to the preservation of $\pm$ matching eigenvalue pairs, only twelve eigenvalues were required to be computed by the rational \textsc{Even-IRA} algorithm. Crosses indicate the two shifts that have been used. Right: 14 eigenvalues (in the upper half plane) found by the rational \textsc{Even-IRA} algorithm (red stars) for $\textnormal{rev} \, P(\lambda)$ with $P(\lambda) = \lambda^2 M + \lambda G + K$ for a gyroscopic system (cf. \cite[Sec.\,4.2]{BenFaStoll}). Blue circles correspond to eigenvalues of $\textnormal{rev} \, P(\lambda)$ computed via \texttt{polyeig}. In contrast to \texttt{polyeig}, the rational \textsc{Even-IRA} algorithm recognizes the fact that the eigenvalues are all located on the imaginary axis.}\label{fig:Plot_Butterfly}
\end{figure}

Our first example is taken from \cite{Butterfly}, see also \texttt{butterfly} in \cite{Bet13}. Here, the matrix polynomial $P(\lambda) = \sum_{j=0}^4 P_j \lambda^j$ under consideration is of degree four. The matrix coefficients are build from several Kronecker products as follows: we set $m=10$ and $n=m^2=100.$ Let $N$ denote the $m \times m$ nilpotent Jordan matrix with ones one the first subdiagonal
and define $\tilde{P}_0 = (1/6)(4 I_m + N + N^T), \tilde{P}_1 = N - N^T,$ $\tilde{P}_2 = -(2I_m - N - N^T)$, $\tilde{P}_3 = \tilde{P}_1$ and $\tilde{P}_4 = - \tilde{P}_2.$ Moreover, we set
$$ P_i = c_{i1} I_m \otimes \tilde{P}_i + c_{i2} \tilde{P}_i \otimes I_m$$
with positive constants $c_{ij}$ chosen as $c_{01}=0.6, c_{02}=1.3, c_{11}=1.3, c_{12}=0.1, c_{21}=0.1, c_{22}=1.2, c_{31}=c_{32}=c_{41}=c_{42}=1.0$ (as in \cite{Butterfly}).  Now the matrix polynomial $P(\lambda) = \sum_{j=0}^4 P_j \lambda^j$ has size $100 \times 100$. We intend to find the 12 eigenvalues of largest magnitude. In order to speed up convergence, we choose a new shift during the iteration as explained in Section \ref{sec:shift} (after a restart) if the first nonzero residual (see \eqref{equ:residual}) is not less than $10^{-5}$. An eigenvalue is considered as converged if its corresponding residual becomes less than $10^{-9}$.
With the initial shift $\zeta \in \CC$ chosen as $0.5 + 2i$, the rational \textsc{Even-IRA} finds the eigenvalues displayed in Figure \ref{fig:Plot_Butterfly} (left plot) in 18 iterations. Compared to the computation with the MATLAB function \texttt{polyeig}, we observe an accordance in both the real and imaginary parts of the computed values of at least the first ten decimal places.  In this experiment, the shift was changed once, see Figure \ref{fig:Plot_Butterfly}, so an LU decomposition of $P(\zeta)$ had to be computed twice (see also the discussion subsequent to Remark \ref{rem:overdetermined}).

\begin{remark}
Recall that the spectral transformation of $\mathcal{L}_P(\zeta) = \zeta X + Y$ to $$K(\zeta) = \mathcal{L}_P(\zeta)^{-T}X \mathcal{L}_P(\zeta)^{-1}X$$ preserves $\pm$ matching pairs of eigenvalues $(+\mu, - \mu)$ as both are mapped to the same eigenvalue $\theta = (\mu^2 - \zeta^2)^{-1}$ (recall Section \ref{sec:EvenIRA}). Therefore, each eigenvalue of $K(\zeta)$ has even multiplicity. In exact arithmetic, multiple eigenvalues will not be captured (see \cite[Sec.\,3]{EvenIRA}) by the Arnoldi iteration. However, as round-off may eventually create them, the authors of the \textsc{Even-IRA} algorithm suggest an additional $X$-orthogonalization of the Krylov basis, see \cite[Lem.\,2.3]{EvenIRA}. Requiring that the basis of the underlying Krylov space is $X$-orthogonal (that is, $\langle v_i , Xv_j \rangle = 0$ for all $i \neq j$) will hinder the algorithm to find multiple copies of the same eigenvalue. The $X$-orthogonalization procedure suggested in \cite{EvenIRA} cannot be directly applied here. The reason for this is that the rational \textsc{Even-IRA} algorithm as outlined in Section \ref{sec:RatEvenIRA} handles complex shifts differently in comparison to \cite{EvenIRA}.
\end{remark}

For our second example we chose the model of a rolling tire, see \cite{ElsVoss} or \cite[Sec.\,4.2.2]{BenFaStoll}. Here $P(\lambda) = \lambda^2 M + \lambda G + K$, where $M,G,K $ are of size $2697 \times 2697$. The matrices $M$ and $K$ are symmetric whereas $G$ is skew-symmetric. Moreover, $M$ and $K$ are positive definite which implies that $P(\lambda)$ has eigenvalues exclusively on the imaginary axis (see \cite[Sec.\,1]{MeerTiss}). Those vary in magnitude from about $10^3$ to $5 \cdot 10^5$. Here we intend to find the eigenvalues of smallest magnitude.  To this end, we consider $\textnormal{rev} \, P(\lambda ) = \lambda^2 K + \lambda G + M$ since the eigenvalues of $\textnormal{rev}\, P(\lambda )$ of largest magnitude correspond via their reciprocals to the eigenvalues of $P(\lambda)$ of smallest magnitude. We have applied the rational \textsc{Even-IRA} to $\textnormal{rev}\, P(\lambda )$ with the same parameters as in the previous example, an initial shift of $10^{-2}i$ and the shift strategy from Section \ref{sec:shift} to find 14 eigenvalues of $\textnormal{rev}\, P(\lambda )$ of largest magnitude.
 The eigenvalues of $\textnormal{rev}\, P(\lambda )$ computed with the rational \textsc{Even-IRA} algorithm and \texttt{polyeig} are displayed in Figure \ref{fig:Plot_Butterfly} (right plot). Eight restarts have been performed. The imaginary parts of the eigenvalues found by the MATLAB function \texttt{polyeig} and those values on the imaginary axis found by the rational \textsc{Even-IRA} algorithm coincide to at least ten significant digits. Clearly, in contrast to \texttt{polyeig}, the rational \textsc{Even-IRA} algorithm anticipates the fact that all eigenvalues are located on the imaginary axis.

\begin{remark}
The shift strategy explained in Section \ref{sec:shift} does not perform optimal for finding the eigenvalues of $P(\lambda)$ with smallest magnitude directly, i.e., when $P(\lambda)$ instead of $\textnormal{rev} \, P(\lambda)$ is used. The shift often increases during the algorithms run and tends to find eigenvalues of larger magnitudes. Thus, according to our experiments, the basic shift strategy from Section \ref{sec:shift} is not appropriate in this situation and a more sophisticated strategy has to be used.
\end{remark}

 In conclusion, the overall success of our algorithm depends in large amounts on the chosen shift-strategy. The method described in Section \ref{sec:shift} works well if one is interested in accelerating the convergence. However, if certain areas of the complex plane are to be ``scanned'' for eigenvalues, a more subtle shift-technique is needed. This is not further discussed here.

\section{Conclusions}
\label{sec:conc}

In this work we have presented a method to compute parts of the spectrum of a $T$-even matrix polynomial. We developed our algorithm on the basis of the \textsc{Even-IRA} algorithm from \cite{EvenIRA} which is a method for computing a few eigenvalues of a $T$-even (i.e. symmetric/skew-symmetric) matrix pencil and the ideas developed in \cite{BennEff} on the rational SHIRA algorithm. Given a $T$-even matrix polynomial, we introduced a special linearization $\mathcal{L}_P(\lambda)$ for $P(\lambda)$ to preserve its $T$-even structure. We showed that the specific block-structure and sparsity of $\mathcal{L}_P(\lambda) = \lambda X + Y$ enables us to solve systems  $\mathcal{L}_P(\zeta)x = y$ in an efficient way. We applied this technique to accelerate the computation of matrix-vector-products for the matrix $K(\zeta) = \mathcal{L}_P(\zeta)^{-T}X \mathcal{L}_P(\zeta)^{-1}X$ to build the underlying Krylov space. An eigenvalue $\theta$ of $K(\zeta)$ gives rise to a $\pm$ matching pair of eigenvalues $+ \sqrt{(1/\theta) + \zeta^2}$ and $ -\sqrt{(1/\theta) + \zeta^2}$ of $\mathcal{L}_P(\lambda)$. As suggested in \cite{EvenIRA}, we used this spectral transformation (i.e. the matrix $K(\zeta)$) and the implicitly restarted Krylov-Schur algorithm to find eigenvalues of $K(\zeta)$. Moreover, we modified the \textsc{Even-IRA} algorithm and turned it into a rational method that is able to handle changes of the shift parameter during the iteration.

A question for future work that is naturally related to our algorithm is the e\-xis\-tence of a compact representation of the Krylov basis similar to the one developed in \cite{Cork} for other types of linearizations (which are not of the same form as $\mathcal{L}_P(\lambda)$). This would be an appropriate means to decrease the cost for storing the Krylov basis vectors.

\section*{Acknowledgments}
Work on this manuscript started when all three authors visited the
Courant Institute of New York University. We would like to give a
special thanks to our host Michael Overton who made this research stay
possible!


\begin{thebibliography}{10}

\bibitem{Templates}
{\sc {Bai, Z., Demmel, J., Dongarra, J., Ruhe, A. and Van der Vorst, H.}}, {\em
  {Templates for the Solution of Algebraic Eigenvalue Problems}}, Society for
  Industrial and Applied Mathematics, Philadelphia, 2000.

\bibitem{Bassour20}
{\sc {Bassour, M.}}, {\em {Hamiltonian polynomial eigenvalue problems}},
  Journal of Applied Mathematics and Physics, 8 (2020), pp.~609 -- 619.

\bibitem{BennEff}
{\sc {Benner, P. and Effenberger, C.}}, {\em {A rational SHIRA method for the
  Hamiltonian eigenvalue problem}}, Taiwanese Journal of Mathematics, 14
  (2010), pp.~805 -- 823.

\bibitem{BenFaStoll}
{\sc {Benner, P., Fassbender, H. and Stoll, M.}}, {\em {Solving large-scale
  quadratic eigenvalue problems with Hamiltonian eigenstructure using a
  structure-preserving Krylov subspace method}}, Electronic Transactions on
  Numerical Analysis, 29 (2008), pp.~212--229.

\bibitem{BenFS11}
{\sc {Benner, P., Fassbender, H. and Stoll, M.}}, {\em { A Hamiltonian
  Krylov–-Schur-type method based on the symplectic Lanczos process}}, Linear
  Algebra and its Applications, 435 (2011), pp.~578 -- 600.

\bibitem{BerlGuett15}
{\sc {Berljafa, M. and G{\"u}ttel, S.}}, {\em {Generalized rational Krylov
  decompositions with an application to rational approximation}}, SIAM Journal
  on Matrix Analysis and Applications, 36 (2015), pp.~894--916.

\bibitem{Bet13}
{\sc {Betcke, T., Higham, N. J., Mehrmann, V. and Tisseur, F.}}, {\em {NLEVP: A
  collection of nonlinear eigenvalue problems}}, ACM Transactions on
  Mathematical Software, 39 (2011).

\bibitem{Demm97}
{\sc {Demmel, J. W.}}, {\em {Applied Numerical Linear Algebra}}, Society for
  Industrial and Applied Mathematics, Philadelphia, 1997.

\bibitem{Dop18}
{\sc {Dopico, F. M., Lawrence, P. W., P{\'e}rez, J. and Van Dooren, P.}}, {\em
  {Block Kronecker linearizations of matrix polynomials and their backward
  errors}}, Numerische Mathematik, 140 (2018), pp.~373--426.

\bibitem{ElsVoss}
{\sc {Elssel, K. and Voss, H.}}, {\em {Reducing huge gyroscopic eigenproblems
  by automated multi-level substructuring}}, Archive of Applied Mechanics, 76
  (2006), pp.~171 -- 179.

\bibitem{FaSa17}
{\sc {Fassbender, H. and Saltenberger, P.}}, {\em {On vector spaces of
  linearizations for matrix polynomials in orthogonal bases}}, Linear Algebra
  and its Applications, 525 (2017), pp.~59 -- 83.

\bibitem{FaSa18}
{\sc {Fassbender, H. and Saltenberger, P.}}, {\em {On a modification of the
  EVEN-IRA algorithm for the solution of T-even polynomial eigenvalue
  problems}}, Proceedings in Applied Mathematics and Mechanics (PAMM),  (2018).

\bibitem{HighMack06}
{\sc {Higham, N. J., Mackey, D. S. and Tisseur, F.}}, {\em {The conditioning of
  linearizations of matrix polynomials}}, SIAM Journal on Matrix Analysis and
  Applications, 28 (2006), pp.~1005 -- 1028.

\bibitem{HighMack06-Sym}
{\sc {Higham, N. J., Mackey, D. S., Mackey, N. and Tisseur, F.}}, {\em
  {Symmetric linearizations for matrix polynomials}}, SIAM Journal on Matrix
  Analysis and Applications, 29 (2006), pp.~143 -- 159.

\bibitem{Reorder}
{\sc {Kagstrom, B.}}, {\em {A direct method for reordering eigenvalues in the
  generalized real Schur form of a regular matrix pair $(A,B)$.}}, In M. S.
  Moonen, G. H. Golub, and B. L. R. De Moor, editors, \emph{Linear Algebra for
  Large Scale and Real-Time Applications}, Kluwer Academic Publishers,
  Amsterdam,  (1993), pp.~195--218.

\bibitem{Mack_Diss}
{\sc {Mackey, D. S.}}, {\em {Structured linearizations for matrix
  polynomials}}, MIMS EPrint 2006.68, Mancheter Institute for Mathematical
  Sciences, Manchester, 2006.

\bibitem{MackPer16}
{\sc {Mackey, D. S. and Perovic, V.}}, {\em {Linearizations of matrix
  polynomials in Bernstein bases}}, Linear Algebra and its Applications, 501
  (2016), pp.~162 -- 197.

\bibitem{MackPer18}
{\sc {Mackey, D. S. and Perovic, V.}}, {\em {Linearizations of matrix
  polynomials in Newton bases}}, Linear Algebra and its Applications, 556
  (2018), pp.~1 -- 45.

\bibitem{MackVibr}
{\sc {Mackey, D. S., Mackey, N., Mehl, C. and Mehrmann, V.}}, {\em {Structured
  polynomial eigenvalue problems: good vibrations from good linearizations}},
  SIAM Journal on Matrix Analysis and Applications, 28 (2006), pp.~1029--1051.

\bibitem{MMMM06}
{\sc {Mackey, D. S., Mackey, N., Mehl, C. and Mehrmann, V.}}, {\em {Vector
  spaces of linearizations for matrix polynomials}}, SIAM Journal on Matrix
  Analysis and Applications, 28 (2006), pp.~971 -- 1004.

\bibitem{MeerTiss}
{\sc {Meerbergen, K. and Tisseur, F.}}, {\em {The quadratic eigenvalue
  problem}}, SIAM Review, 43 (2001), pp.~235--286.

\bibitem{MehrWat01}
{\sc {Mehrmann, V. and Watkins, D.}}, {\em {Structure-preserving methods for
  computing eigenpairs of large sparse skew-Hamiltoninan/Hamiltonian pencils}},
  SIAM Journal on Scientific Computing, 22 (2001), pp.~1905--1925.

\bibitem{Butterfly}
{\sc {Mehrmann, V. and Watkins, D.}}, {\em {Polynomial eigenvalue problems with
  Hamiltonian structure}}, Electronic Transactions on Numerical Analysis, 13
  (2002), pp.~106 -- 118.

\bibitem{EvenIRA}
{\sc {Mehrmann, V., Schr{\"o}der, C. and Simoncini, V.}}, {\em {An
  implicitly-restarted Krylov subspace method for real symmetric/skew-symmetric
  eigenproblems}}, Linear Algebra and its Applications, 436 (2012),
  pp.~4070--4087.

\bibitem{MolStew73}
{\sc {Moler, C. B. and Stewart, G. W.}}, {\em {An algorithm for generalized
  matrix eigenvalue problems}}, SIAM Journal on Numerical Analysis, 10 (1973),
  pp.~241--256.

\bibitem{Ruhe84}
{\sc {Ruhe, A.}}, {\em {Rational Krylov sequence methods for eigenvalue
  computation}}, Linear Algebra and its Applications, 58 (1984), pp.~391--405.

\bibitem{Ruhe98}
{\sc {Ruhe, A.}}, {\em {Rational Krylov: a practical algorithm for large sparse
  nonsymmetric matrix pencils}}, SIAM Journal on Scientific Computing, 19
  (1998), pp.~1535--1551.

\bibitem{PS_Diss}
{\sc {Saltenberger, P.}}, {\em {On different concepts for the linearization of
  matrix polynomials and canonical decompositions of structured matrices with
  respect to indefinite sesquilinear forms}}, Logos Verlag, Berlin, 2019.

\bibitem{KrylovSchur}
{\sc {Stewart, G. W.}}, {\em {A Krylov–-Schur algorithm for large
  eigenproblems}}, SIAM Journal on Matrix Analysis and Applications, 23 (2002),
  pp.~601--614.

\bibitem{Cork}
{\sc {Van Beeumen, R., Meerbergen, K. and Michiels, W.}}, {\em {Compact
  rational Krylov methods for nonlinear eigenvalue problems}}, SIAM Journal on
  Matrix Analysis and Applications, 36 (2015), pp.~820 -- 838.

\bibitem{Wat04}
{\sc {Watkins, D.}}, {\em {On Hamiltonian and symplectic Lanczos processes}},
  Linear Algebra and its Applications, 385 (2004), pp.~23--45.

\end{thebibliography}
\end{document}